# VOTER MODELS WITH HETEROZYGOSITY SELECTION

By Anja Sturm[1] and Jan Swart[2]

*University of Delaware and ÚTIA Prague*

This paper studies variations of the usual voter model that favor types that are locally less common. Such models are dual to certain systems of branching annihilating random walks that are parity preserving. For both the voter models and their dual branching annihilating systems we determine all homogeneous invariant laws, and we study convergence to these laws started from other initial laws.

**1. Introduction and main results.**

1.1. *Voter models with heterozygosity selection.* This paper studies variations of the usual voter model that favor types that are locally less common. These systems can be used to model the distribution of two types of organisms (two similar species or merely different genetic variants of the same species) that occupy overlapping ecological niches, and therefore compete with each other for resources. If both types are equally fit, but their ecological niches are not completely identical, then individuals belonging to the type that is locally less common have an advantage, since they can use resources that are not used by most of their neighbors. This effect is called *negative frequency dependent selection* or (positive) *heterozygosity selection*. (Here, following a common practice in population biology, the word heterozygosity refers to the degree of genetic variation within a population as a whole, rather than the variation between homologous chromosomes in a diploid organism.)

Our processes of interest are Markov processes $X = (X_t)_{t \geq 0}$ with state space $\{0,1\}^{\mathbb{Z}^d}$. We denote a typical element of $\{0,1\}^{\mathbb{Z}^d}$ by $x = (x(i))_{i \in \mathbb{Z}^d}$,

Received January 2007; revised April 2007.
[1]Supported by GAČR Grants 201/06/1323 and 201/07/0237.
[2]Supported in part by UDRF Grant 06000596.

where $x(i) \in \{0, 1\}$ is interpreted as the *type* of the organism at the site $i$. Borrowing terminology from physics, we sometimes also call $x(i)$ the *spin* at $i$. We say that a $\{0,1\}^{\mathbb{Z}^d}$-valued Markov process $X$ is *spin-flip symmetric* if its dynamics are symmetric under a simultaneous flip of all spins, that is, the transition $x \mapsto x'$ happens at the same rate as the transition $(1-x) \mapsto (1-x')$.

More specifically, we are interested in the following models.

DEFINITION 1 (*Neutral Neuhauser–Pacala model*). The *neutral Neuhauser–Pacala model* is the spin-flip symmetric Markov process in $\{0,1\}^{\mathbb{Z}^d}$ such that if the state of the process is $x$, then $x(i)$ flips from 0 to 1 with rate $f_1(f_0 + \alpha f_1)$ (and likewise for jumps from 1 to 0, by spin-flip symmetry), where

$$(1.1) \qquad f_\tau := \frac{1}{|\mathcal{N}_i|} \sum_{j \in \mathcal{N}_i} 1_{\{x(j) = \tau\}} \qquad (\tau = 0, 1)$$

denotes the local frequency of type $\tau$ in the block of $(2R+1)^d - 1$ sites centered around $i$, not containing $i$ itself, given by

$$(1.2) \qquad \mathcal{N}_i := \{j = (j_1, \ldots, j_d) \in \mathbb{Z}^d : 0 < |j_k - i_k| \leq R \ \forall k\},$$

and $0 \leq \alpha \leq 1$, $d, R \geq 1$ are parameters such that $\min\{d, R\} \geq 2$.

The neutral Neuhauser–Pacala model is a special case of the model introduced in [21], when their parameters satisfy $\lambda = 1$ and $\alpha_{01} = \alpha_{10} =: \alpha$ (called the "symmetric case" there). For $\alpha = 1$, this is a usual (range $R$) voter model, while for $\alpha < 1$, locally rare types have an advantage. Neuhauser and Pacala interpret the rate $f_1(f_0 + \alpha f_1)$ as follows: at each site $i$, an organism of type 0 dies with rate $f_0 + \alpha f_1$ due to competition with its neighbors, and is immediately replaced by an organism of a random type chosen from $\mathcal{N}_i$. If $\alpha < 1$, then the *interspecific competition* is smaller than the *intraspecific competition*, hence, locally rare types die less frequently.

Our next model of interest is another nonlinear voter model.

DEFINITION 2 (*Affine voter model*). The *affine voter model* is the spin-flip symmetric Markov process in $\{0,1\}^{\mathbb{Z}^d}$ such that if the state of the process is $x$, then $x(i)$ flips from 0 to 1 with rate $\alpha f_1 + (1-\alpha) 1_{\{f_1 > 0\}}$, where $f_0, f_1$ are defined as for the previous model. Here $0 \leq \alpha \leq 1$ and $d, R \geq 1$, $\min\{d, R\} \geq 2$.

Again, for $\alpha = 1$, this is a usual (range $R$) voter model, while for $\alpha < 1$, locally rare types have an advantage. For $\alpha = 0$, the affine voter model is a threshold voter model [7, 14, 18]. In this and the previous definition, we have excluded the case $d = 1 = R$ since this model has special behavior (see Section 2.1).

The third and final model we will consider is a one-dimensional model.

DEFINITION 3 (*Rebellious voter model*). The *rebellious voter model* is the spin-flip symmetric Markov process in $\{0,1\}^{\mathbb{Z}}$ such that if the state of the process is $x$, then $x(i)$ flips with rate

$$
\begin{aligned}
&\alpha(1_{\{x(i-1)\neq x(i)\}} + 1_{\{x(i)\neq x(i+1)\}}) \\
&\quad + (1-\alpha)(1_{\{x(i-2)\neq x(i-1)\}} + 1_{\{x(i+1)\neq x(i+2)\}}).
\end{aligned}
\tag{1.3}
$$

For $\alpha = 1$, this is a nearest-neighbor one-dimensional voter model, but for $\alpha < 1$, locally rare types have an advantage. To see this, note that if $x(i+1) \neq x(i+2)$ and the spin at $i$ flips, then $x(i)$ was previously of the most common type in the set $\{i, i+1, i+2\}$.

In all of the models above, we call $\alpha$ the (interspecific) *competition parameter*. As usual with voter models, the main interest in these models lies in the phase transition between coexistence and noncoexistence. It is believed, and has been proved in special cases, that coexistence occurs for sufficiently small competition parameters $\alpha$ or in high dimensions, while noncoexistence occurs in low dimensions and for large competition parameters.

We note that Blath, Etheridge and Meredith [1] have studied systems of interacting Wright–Fisher diffusions with heterozygosity selection, that is, the stepping stone analogue of the voter models discussed in the present paper.

1.2. *Duality with parity preserving processes.* The neutral Neuhauser–Pacala model and the affine and rebellious voter models are cancellative spin systems. Here, a Markov process $X = (X_t)_{t \geq 0}$ with state space $\{0,1\}^{\mathbb{Z}^d}$ is a *cancellative spin system* if its dynamics are of the following special form. For each finite set $A \subset \mathbb{Z}^d \times \mathbb{Z}^d$, there is a rate $a(A) \geq 0$, such that with this rate, the process jumps from the state $x$ to $x + Ax \mod(2)$, where

$$
Ax(i) := \sum_{j\,:\,(i,j) \in A} x(j) \mod(2) \quad (i \in \mathbb{Z}^d).
\tag{1.4}
$$

(In [13], a somewhat more general class of cancellative spin systems is considered, where also spontaneous flips are allowed.)

The neutral Neuhauser–Pacala model can be cast into the general form of a cancellative spin system, with rates $a(A)$ given by

$$
\begin{aligned}
a(\{i\} \times \{i,j\}) &= \frac{\alpha}{|\mathcal{N}_i|} \qquad (j \in \mathcal{N}_i), \\
a(\{i\} \times \{j,k\}) &= \frac{1-\alpha}{|\mathcal{N}_i|^2} \qquad (j,k \in \mathcal{N}_i,\ j \neq k)
\end{aligned}
\tag{1.5}
$$

and $a(A) = 0$ in all other cases. Likewise, the affine voter model can be formulated as a cancellative spin system, with rates $a(A)$ given by

$$
\begin{aligned}
a(\{i\} \times \{i,j\}) &= \alpha |\mathcal{N}_i|^{-1} & (j \in \mathcal{N}_i), \\
a(\{i\} \times \Delta) &= (1-\alpha) 2^{-|\mathcal{N}_i|+1} & (\Delta \subset \mathcal{N}_i \cup \{i\}, |\Delta| \text{ even}),
\end{aligned}
\tag{1.6}
$$

while in case of the rebellious voter model, the rates are given by

$$
\begin{aligned}
a(\{i\} \times \{i-1,i\}) &= a(\{i\} \times \{i,i+1\}) = \alpha, \\
a(\{i\} \times \{i-2,i-1\}) &= a(\{i\} \times \{i+1,i+2\}) = 1-\alpha,
\end{aligned}
\qquad (i \in \mathbb{Z})
\tag{1.7}
$$

and $a(A) = 0$ in all other cases.

Under a suitable summability assumption on the rates, for each cancellative spin system $X$ there exists a unique cancellative spin system $Y$ such that $X$ and $Y$ are dual to each other in the sense that

$$
\mathbb{P}[|X_t Y_0| \text{ is odd}] = \mathbb{P}[|X_0 Y_t| \text{ is odd}] \qquad (t \geq 0),
\tag{1.8}
$$

whenever $X$ and $Y$ are independent (with arbitrary initial laws), and either $|X_0|$ or $|Y_0|$ is finite (see [13]). Here, for any $x, y \in \{0,1\}^{\mathbb{Z}}$, we write $|x| := \sum_i x(i)$ and $(xy)(i) := x(i)y(i)$. If $X$ is defined by rates $a_X(A)$, then $Y$ is defined by rates $a_Y(A)$ given by $a_Y(A) = a_X(A^{\mathrm{T}})$, where $A^{\mathrm{T}} := \{(j,i) : (i,j) \in A\}$. We note that since the functions $y \mapsto 1_{\{|xy| \text{ is odd}\}}$ with $|x| < \infty$ are distribution determining, (1.8) determines the transition laws of $Y$ uniquely.

For example, if $X$ is the rebellious voter model, then the dynamics of $Y$ have the following description. We interpret the sites $i$ for which $Y_t(i) = 1$ as being occupied by a particle at time $t$. Then each particle jumps with rate $\alpha$ one step to the left, and with the same rate to the right, with the rule that if the site the particle lands on is already occupied, the two particles annihilate. Moreover, with rate $1 - \alpha$, each particle gives birth to two new particles located on its nearest and next-nearest site to the left, and with the same rate on the right, again annihilating with any particles that may already be present on these sites. We call $Y$ the *asymmetric double branching annihilation random walk* (ADBARW). Also in case of the neutral Neuhauser–Pacala and affine voter models, one may check that the dual model is a *system of branching annihilation random walks*, with branching rate proportional to $1 - \alpha$. Since the number of particles that are born is always even, these systems preserve parity.

More generally, if $X$ is a cancellative spin system and $Y$ is its dual, let us say that $Y$ is *parity preserving* if the process started in any finite initial state satisfies $(-1)^{|Y_t|} = (-1)^{|Y_0|}$ a.s. It is not hard to see that the following statements are equivalent: (i) $X$ is spin-flip symmetric, (ii) $Y$ is parity preserving, (iii) $a_X(A) = 0$ unless $|\{j : (i,j) \in A\}|$ is even for all $i \in \mathbb{Z}^d$.

We next define the concepts of *coexistence*, *persistence*, *survival* and *stability*. Below, $\underline{0}$ and $\underline{1}$ denote the configurations in $\{0,1\}^{\mathbb{Z}^d}$ which are identically zero and one, respectively.

DEFINITION 4 (*Coexistence*). We say that a probability law $\mu$ on $\{0,1\}^{\mathbb{Z}^d}$ is *coexisting* if $\mu(\{\underline{0},\underline{1}\}) = 0$. We say that a spin-flip symmetric cancellative spin system $X$ exhibits *coexistence* if there exists a coexisting invariant law for $X$.

DEFINITION 5 (*Persistence*). We say that a probability law $\mu$ on $\{0,1\}^{\mathbb{Z}^d}$ is *nonzero* if $\mu(\{\underline{0}\}) = 0$. We say that a parity preserving cancellative spin system $Y$ exhibits *persistence* if there exists a nonzero invariant law for $Y$.

Below, $X^x$ and $Y^y$ denote the processes $X$ and $Y$ started in $X_0^x = x$ and $Y_0^y = y$, respectively.

DEFINITION 6 (*Survival*). We say that a spin-flip symmetric cancellative spin system $X$ *survives* if

$$(1.9) \qquad \mathbb{P}[X_t^x \neq \underline{0} \ \forall t \geq 0] > 0 \qquad \text{for some } |x| < \infty.$$

We say that a parity preserving cancellative spin system $Y$ *survives* if

$$(1.10) \qquad \mathbb{P}[Y_t^y \neq \underline{0} \ \forall t \geq 0] > 0 \qquad \text{for some } |y| < \infty, |y| \text{ even}.$$

Note that because of parity preservation, the left-hand side of (1.10) is always one if $|y|$ is odd. In view of this, it will usually be clear what we mean when we say that a process survives. Our definitions are not entirely unambiguous, however, since it is possible for a cancellative spin system to be both spin-flip symmetric and parity preserving. When there is danger of confusion, we will say that the *even process survives* if we mean survival in the sense of (1.10).

For a parity preserving cancellative spin system $Y$ started in an odd initial state, we define

$$(1.11) \qquad \hat{Y}_t(i) := Y_t(l(t) + i) \qquad (i \in \mathbb{N}^d)$$
$$\text{where } l(t) := \inf\{i \in \mathbb{Z}^d : Y_t(i) = 1\}.$$

Here, the infimum is defined componentwise. Note that $(\hat{Y}_t)_{t \geq 0}$ is a Markov process with state space $\{y \in \{0,1\}^{\mathbb{N}^d} : |y| \text{ is finite and odd}, \inf\{i : y(i) = 1\} = \underline{0}\}$. We call $\hat{Y}$ the process $Y$ *viewed from its lower left corner*.

DEFINITION 7 (*Stability*). We say that a parity preserving cancellative spin system $Y$ is *stable* if the state with one particle at the origin is positively recurrent for the Markov process $\hat{Y}$.

Using duality, one can prove that each cancellative spin system has a special invariant law, which is the limit law of the process started in product law with intensity $1/2$ [13]. Thus,

$$
\begin{aligned}
\mathbb{P}[X_t^{1/2} \in \cdot] &\underset{t \to \infty}{\Longrightarrow} \mathbb{P}[X_\infty^{1/2} \in \cdot] =: \nu_X^{1/2}, \\
\mathbb{P}[Y_t^{1/2} \in \cdot] &\underset{t \to \infty}{\Longrightarrow} \mathbb{P}[Y_\infty^{1/2} \in \cdot] =: \nu_Y^{1/2},
\end{aligned}
\tag{1.12}
$$

where $\nu_X^{1/2}$ and $\nu_Y^{1/2}$ are invariant laws for the processes $X$ and $Y$, respectively. Because of certain analogies with additive spin systems, we call $\nu_X^{1/2}$ and $\nu_Y^{1/2}$ the *odd upper invariant laws* of $X$ and $Y$, respectively.

The next lemma is a simple consequence of duality.

LEMMA 1 (Invariant laws and survival). *Let $X$ be a spin-flip symmetric cancellative spin system and let $Y$ be its dual parity preserving cancellative spin system. Then:*

(a) *The following statements are equivalent:* (i) $X$ *exhibits coexistence*, (ii) $\nu_X^{1/2}$ *is not concentrated on* $\{\underline{0}, \underline{1}\}$, (iii) $Y$ *survives*.

(b) *The following statements are equivalent:* (i) $Y$ *persists*, (ii) $\nu_Y^{1/2}$ *is not concentrated on* $\underline{0}$, (iii) $X$ *survives*.

Apart from the relations in Lemma 1, for our models of interest, one readily conjectures a number of other relations between the concepts we have just introduced. While being supported by numerical simulations, these conjectures appear to be hard to prove in general. Thus, we formulate as open questions:

Q1. Is coexistence of $X$ equivalent to survival of $X$?
Q2. Does stability of $Y$ imply extinction of $Y$?
Q3. Does coexistence for $\alpha$ imply coexistence for all $\alpha' < \alpha$?

The rebellious voter model has a special property, explained in Section 2.1, that allows us to answer Question Q1 positively for this model.

LEMMA 2 (Survival and coexistence). *The rebellious voter model survives if and only if it exhibits coexistence.*

1.3. *Results.* We say that a probability law $\mu$ on $\{0,1\}^{\mathbb{Z}^d}$ is *homogeneous* if $\mu$ is translation invariant. The odd upper invariant laws $\nu_X^{1/2}$ and $\nu_Y^{1/2}$ are examples of homogeneous invariant laws, and so are the delta-measures $\delta_{\underline{0}}$ and (in case of $X$) $\delta_{\underline{1}}$. We will see that, under weak additional assumptions, there are no others.

For additive spin systems, duality may be employed to show that, under certain conditions, the upper invariant law is the only nonzero homogeneous invariant law, and the long-time limit law for any system started from a nonzero homogeneous initial law; see [15] and [17], Theorem III.5.18. With certain complications, these techniques can be adapted to cancellative spin systems and their odd upper invariant laws. This idea has been successfully applied in [4] to certain annihilating branching processes. As our next theorem demonstrates, it can be made to work for our models as well.

THEOREM 3 (Homogeneous invariant laws). *Let $X$ be either the neutral Neuhauser–Pacala model, the affine voter model or the rebellious voter model and let $Y$ be its dual. Then:*

(a) *If $\alpha < 1$ and $Y$ survives, then $\nu_X^{1/2}$ is the unique homogeneous coexisting invariant law of $X$. If, moreover, $\alpha > 0$ and $Y$ is not stable, then the process $X$ started in any homogeneous coexisting initial law satisfies*

$$\mathbb{P}[X_t \in \cdot] \underset{t \to \infty}{\Longrightarrow} \nu_X^{1/2}. \tag{1.13}$$

(b) *If $\alpha < 1$ and $X$ survives, then $\nu_Y^{1/2}$ is the unique homogeneous nonzero invariant law of $Y$. If, moreover, $\alpha > 0$ and $d \geq 2$, then the process $Y$ started in any homogeneous nonzero initial law satisfies*

$$\mathbb{P}[Y_t \in \cdot] \underset{t \to \infty}{\Longrightarrow} \nu_Y^{1/2}. \tag{1.14}$$

Theorem 3 does not tell us anything about the values of $\alpha$ for which coexistence occurs. For the affine voter model, coexistence for $\alpha = 0$ has been proved in [7, 18]. It seems likely that, using comparison with oriented percolation, this result can be extended to small positive $\alpha$. Comparison with oriented percolation was used in [21], where coexistence for the neutral Neuhauser–Pacala model for $\alpha$ sufficiently close to zero is proved, and in [9] where (among other things) coexistence is proved for the neutral Neuhauser–Pacala model in dimensions $d \geq 3$ for $\alpha$ sufficiently close to *one*.

Again using comparison with oriented percolation, we will prove that, for small $\alpha$, the rebellious voter model coexists. In fact, we prove considerably more.

THEOREM 4 (Complete convergence). *Let $X$ be the rebellious voter model and let $Y$ be its dual. Then there exists an $\alpha' > 0$, such that for all $\alpha \in [0, \alpha')$:*

(a) *The process $X$ exhibits coexistence and survival. The process $Y$ exhibits persistence, survival and is not stable.*

(b) *The process $X$ started in an arbitrary initial law satisfies*

$$\text{(1.15)} \qquad \mathbb{P}[X_t \in \cdot] \underset{t \to \infty}{\Longrightarrow} \rho_0 \delta_{\underline{0}} + \rho_1 \delta_{\underline{1}} + (1 - \rho_0 - \rho_1) \nu_X^{1/2},$$

*where* $\rho_q := \mathbb{P}[X_t = \underline{q} \text{ for some } t \geq 0]$ $(q = 0, 1)$.

In analogy with similar terminology for the contact process, we call the result in part (b) *complete convergence*. Complete convergence for the threshold voter model (i.e., the affine voter model with $\alpha = 0$) was proved in [14].

## 2. Methods and discussion.

2.1. *Interfaces.* For one-dimensional models, there is, apart from duality, a useful additional tool available. If $X$ is a spin-flip symmetric cancellative spin system on $\mathbb{Z}$, then setting

$$\text{(2.1)} \qquad Y_t(i) := 1_{\{X_t(i) \neq X_t(i+1)\}} \qquad (t \geq 0, i \in \mathbb{Z})$$

defines a Markov process $Y = (Y_t)_{t \geq 0}$ in $\{0,1\}^{\mathbb{Z}}$ that we call the *interface model* associated with $X$. Under a suitable summability assumption on the rates, $Y$ is a parity preserving cancellative spin system.

By looking at interfaces, we can explain our interest in the rebellious voter model. Moreover, we can explain why we have excluded the case $d = 1 = R$ from Definitions 1 and 2. Consider the following spin-flip symmetric models.

DEFINITION 8 (*One-dimensional models*). For any $0 \leq \alpha \leq 1$, the *disagreement voter model* is the cancellative spin system on $\mathbb{Z}$ with rates $a(A)$ given by

$$\text{(2.2)} \qquad \begin{aligned} a(\{i\} \times \{i-1, i\}) = a(\{i\} \times \{i, i+1\}) &= \alpha, \\ a(\{i\} \times \{i-1, i+1\}) &= 1 - \alpha, \end{aligned} \qquad (i \in \mathbb{Z}),$$

and the *swapping voter model* is given by the rates

$$\text{(2.3)} \qquad \begin{aligned} a(\{i\} \times \{i-1, i\}) = a(\{i\} \times \{i, i+1\}) &= \alpha, \\ a(\{i, i+1\} \times \{i, i+1\}) &= 1 - \alpha, \end{aligned} \qquad (i \in \mathbb{Z}).$$

In each case, it is understood $a(A) = 0$ for all $A$ other than those mentioned.

If in the definition of the neutral Neuhauser–Pacala model with competition parameter $\alpha_{\text{NP}}$ one would set $d = 1 = R$, then up to a trivial redefinition of the speed of time, one would obtain a disagreement voter model with parameter $\alpha = 2\alpha_{\text{NP}}/(1 + \alpha_{\text{NP}})$. Likewise, setting $d = 1 = R$ in the definition of the affine voter model with competition parameter $\alpha_{\text{AV}}$ yields, up to a change of the speed of time, a disagreement voter model with parameter $\alpha = 1/(2 - \alpha_{\mathrm{AV}})$.

To explain the special properties of the disagreement and rebellious voter models, we look at the way these models are related to other models through duality and interface relations. These relations are summarized in Figure 1. Recall that the dual of the rebellious voter model is the ADBARW. The dual of the disagreement voter model has been called the *double branching annihilating random walk* (DBARW) in [24]. The swapping voter model has a mixture of voter model and exclusion process dynamics. We call its dual the *swapping annihilating random walk* (SARW).

The swapping annihilating random walk (SARW) has the special property that the number of particles cannot increase. As a result, the behavior of the disagreement and swapping voter models and their duals is largely known. For each $\alpha > 0$, the disagreement and swapping voter models exhibit extinction and noncoexistence. The DBARW and the SARW get extinct and are not persistent. These facts follow from [24], Theorem 8 and [21], Theorem 2(b). The disagreement voter model with $\alpha = 1/2$ is a threshold voter model and has been studied earlier in [7]. It is trivial that the SARW is stable. We will prove stability for the DBARW elsewhere [23].

The rebellious voter model has the special property that its interface model and dual coincide, which is one of our main reasons for introducing it. This property is very helpful in the proof of Theorem 4. Moreover, it allows us to prove Lemma 2:

PROOF OF LEMMA 2. It is easy to see that $X$ survives if and only if its (even) interface model $Y$ survives. Since $Y$ is also the dual of $X$, by Lemma 1, $X$ exhibits coexistence if and only if $Y$ survives. □

For the rebellious voter model, numerical simulations suggest the existence of a critical value $\alpha_{\mathrm{c}} \approx 0.5$ such that, for $\alpha < \alpha_{\mathrm{c}}$, one has survival, coexistence and instability of the interface model, while, for $\alpha > \alpha_{\mathrm{c}}$, one has extinction, noncoexistence and stability of the interface model (see Figure 2). We conjecture qualitatively similar behavior (but with a different critical point) for the one-dimensional neutral Neuhauser–Pacala and affine voter models. Note that in this respect, these models are quite different from the disagreement voter model, which gets extinct for all $\alpha > 0$.

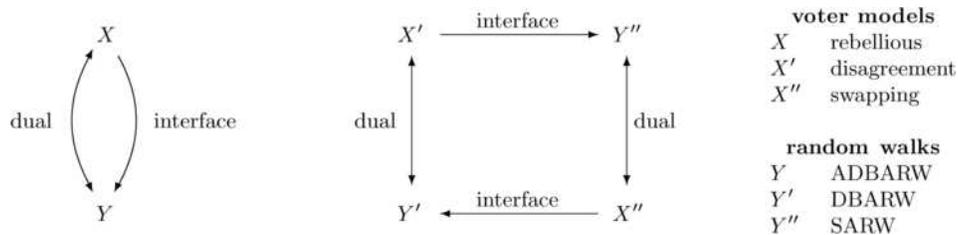

FIG. 1. *Relations between models.*

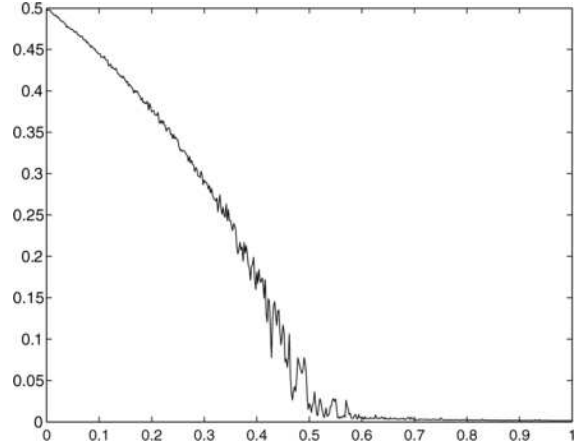

Fig. 2. *Equilibrium density of the ADBARW as a function of $\alpha$. The data were obtained starting with one particle on an interval of 700 sites with periodic boundary conditions, gradually lowering $\alpha$ from one to zero during a time interval of length 300,000.*

2.2. *Homogeneous laws.* In this section we discuss the methods used to prove Theorem 3. Using duality, it is not hard to see (see Section 3.1) that the odd upper invariant law of $X$, defined in (1.12), is uniquely determined by

(2.4) $\qquad \mathbb{P}[|X_\infty^{1/2} y| \text{ is odd}] = \frac{1}{2}\mathbb{P}[Y_s^y \neq \underline{0} \ \forall s \geq 0] \qquad (|y| < \infty).$

In order to prove the convergence in (1.13), it therefore suffices to show that

(2.5) $\qquad \lim_{t \to \infty} \mathbb{P}[|X_t y| \text{ is odd}] = \frac{1}{2}\mathbb{P}[Y_s^y \neq \underline{0} \ \forall s \geq 0] \qquad (|y| < \infty).$

Using duality (1.8), we can rewrite the left-hand side as

(2.6) $\qquad \lim_{t \to \infty} \mathbb{P}[|X_{t_0} Y_{t-t_0}^y | \text{ is odd}],$

where $t_0 > 0$ is fixed. Conditioning on the event of survival, we need to show that

(2.7) $\qquad \lim_{t \to \infty} \mathbb{P}[|X_{t_0} Y_{t-t_0}^y | \text{ is odd} \,|\, Y_s^y \neq \underline{0} \ \forall s \geq 0] = \frac{1}{2}.$

It turns out that we can show this, provided that we can show that

(2.8) $\qquad \lim_{t \to \infty} \mathbb{P}[0 < |Y_t^y| < N] = 0 \qquad (N \geq 1),$

that is, $Y$ exhibits a form of *extinction versus unbounded growth*. To prove this, we need to assume that $Y$ is not stable. The proof uses induction on $N$ and is quite subtle. Originally, we only knew how to prove convergence in Césaro mean, until we saw the paper of Handjani [14] where usual convergence as $t \to \infty$ is proved in the context of threshold voter models.

The proof of the convergence in (1.14) follows similar lines. This time, instead of (2.7), we need to prove

$$(2.9) \qquad \lim_{t \to \infty} \mathbb{P}[|X^x_{t-t_0} Y_{t_0}| \text{ is odd} \mid X^x_s \neq \underline{0} \; \forall s \geq 0] = \tfrac{1}{2}.$$

At first, it might seem that this is true provided that, in analogy with (2.8), one has that $\mathbb{P}[0 < |X^x_t| < N] \to 0$ for all $N \geq 1$. Indeed, the reference [22] contains a claim of this sort, but as we explain in Section 3.6 below, this is not correct. Indeed, in order for the probability on the left-hand side of (2.9) to be close to $\tfrac{1}{2}$, we need many events that could affect the parity of $|X^x_{t-t_0} Y_{t_0}|$. This means that there must be many sites $i, j, k$, close to each other, such that $Y_0(i) = 1$, while $X^x_{t-t_0}(j) \neq X^x_{t-t_0}(k)$. Indeed, it suffices to verify that, conditional on survival, the quantity

$$(2.10) \qquad |\nabla X^x_t| := |\{(i,j) \in \mathbb{Z}^d \times \mathbb{Z}^d : |i-j| = 1, X^x_t(i) \neq X^x_t(j)\}|$$

tends to infinity as $t \to \infty$. In dimensions $d \geq 2$, we can verify this, but in dimension $d = 1$, we run into the difficulty that it is hard to rule out the scenario that, at certain large times, $X^x_t$ consists of just one large interval of ones, in which case $|\nabla X^x_t| = 2$. (Indeed, to prove a result in $d = 1$, one would probably need to assume that the interface model associated with $X$ is not stable. We have not carried out this approach.)

2.3. *Comparison with oriented percolation.* Theorem 4 is similar to the main result of [14]. In that paper a key technical tool is comparison of threshold voter models with threshold contact processes, which in [18] were shown to survive. Using complete convergence for these threshold contact processes, one can then prove the analogue statement for the threshold contact processes.

Our approach will be similar, except that we will use comparison with oriented percolation. Also, instead of proving a comparison result for the rebellious voter model $X$, we prefer to work with the ADBARW $Y$, which is both the dual and the interface model associated with $X$. We will prove that if $\alpha$ is small, then for each $p < 1$, the process $Y$ viewed on suitable length and time scales dominates an oriented percolation with parameter $p$. This kind of rescaling argument was first used in Bramson and Durrett [3]. In [4, 11], it was shown how the technique may be amended to cover also interacting particle systems that lack monotonicity. Unfortunately, these references are somewhat imprecise when it comes to showing $m$-dependence for the objects one compares with; this is done more carefully in [12]. We will use a somewhat different argument to ensure $m$-dependence than the one used in that reference.

Since the comparison result is of some interest on its own, we formulate it here. We first introduce oriented (site) percolation with percolation parameter $p$. Let $\mathbb{Z}^2_{\text{even}} := \{(x,n) \in \mathbb{Z}^2 : x + n \text{ is even}\}$. Let $\{\omega_z : z \in \mathbb{Z}^2_{\text{even}}\}$ be i.i.d. Bernoulli random variables with $\mathbb{P}[\omega_z = 1] = p$. For $z, z' \in \mathbb{Z}^2_{\text{even}}$, we say that there is an open path from $z$ to $z'$, denoted as $z \to z'$, if there exist $(x_n, n), \ldots, (x_m, m) \in \mathbb{Z}^2_{\text{even}}$ with $|x_k - x_{k-1}| = 1$ and $\omega_{(x_k, k)} = 1$ for all $n < k \leq m$, such that $(x_n, n) = z$ and $(x_m, m) = z'$. By definition, $z \to z$ for all $z \in \mathbb{Z}^2_{\text{even}}$. For given $A \subset \mathbb{Z}_{\text{even}} := \{2n : n \in \mathbb{Z}\}$, we put for $n \geq 0$

$$(2.11) \quad W_n := \{x \in \mathbb{Z} : (x, n) \in \mathbb{Z}^2_{\text{even}}, \exists x' \in A \text{ s.t. } (x', 0) \to (x, n)\}.$$

Then $W = (W_n)_{n \geq 0}$ is a Markov chain, taking values, in turn, in the subsets of $\mathbb{Z}_{\text{even}}$ and $\mathbb{Z}_{\text{odd}}$, started in $W_0 = A$. We call $W$ the *oriented percolation process*.

The comparison entails defining certain "good" events concerning the behavior of the ADBARW in large space-time boxes. Let $L \geq 1$ and $T > 0$. For any $x \in \mathbb{Z}$, put

$$(2.12) \quad \begin{aligned} I_x &:= \{2Lx - L, \ldots, 2Lx + L\} \quad \text{and} \\ I'_x &:= \{2Lx - 4L, \ldots, 2Lx + 4L\}. \end{aligned}$$

Let $Y$ be an ADBARW started in an arbitrary initial state $Y_0$. We define a set of "good" points for $n \geq 0$ by

$$(2.13) \quad \begin{aligned} \chi_n &:= \{x \in \mathbb{Z} : (x, n) \in \mathbb{Z}^2_{\text{even}}, \exists i \in I_x \text{ s.t. } Y_{nT}(i) = 1 \text{ and,} \\ &\quad \text{in case } n \geq 1, \forall (n-1)T < t \leq nT \ \exists i \in I'_x \text{ s.t. } Y_t(i) = 1\}. \end{aligned}$$

With these definitions, our result reads:

THEOREM 5 (Comparison with oriented percolation). *For each $p < 1$, there exists an $\alpha' > 0$ such that for all $\alpha \in [0, \alpha')$ there exist $L \geq 1$ and $T > 0$, such that if $Y$ is an ADBARW with parameter $\alpha$, started in an arbitrary initial state $Y_0$, then the process $(\chi_n)_{n \geq 0}$ defined in (2.13) can be coupled to an oriented percolation process $(W_n)_{n \geq 0}$ with parameter $p$ and initial state $W_0 = \chi_0$, in such a way that $W_n \subset \chi_n$ for all $n \geq 0$.*

2.4. *Discussion and open problems.* Three open problems have already been formulated in Section 1.2. As a fourth problem, we mention the following:

Q4. For which values of the parameters do voter models with heterozygosity selection exhibit noncoexistence?

As mentioned before, there are several results proving coexistence for the sort of models we are considering. Very little is known about noncoexistence, except for the pure voter model in dimensions $d = 1, 2$ and the disagreement voter model, which is somehow special. There is some hope that the methods used in [9] (see also [8]), who prove coexistence for the neutral Neuhauser–Pacala model with $\alpha$ close to one in dimensions $d \geq 3$, can be extended to cover dimension 2 as well. Therefore, noncoexistence can be expected in dimension one only. In fact, it seems that physicists believe that for these and similar models, there is a critical dimension $d_c \approx 4/3$ such that only below $d_c$ there is a nontrivial phase transition between coexistence and noncoexistence [5, 6, 25]. Here, the fractional dimension probably refers to self-similar lattices.

Note that, by duality, proving noncoexistence for large interspecific competition boils down to proving extincton for a parity preserving branching process with a small branching rate. Usually, for interacting particle systems, it is easier to find sufficient conditions for extinction than for survival; this is true for the contact process, and also for the annihilating branching process studied in [4]. The difficulties in our case come from parity preservation, which makes extinction difficult, slow (slower than exponential) and "nonlocal," since it may require particles to come from far away to annihilate each other.

2.5. *Outline.* The rest of the paper is devoted to proofs. Lemma 2 has already been proved in Section 2.1 above. In Section 3 below, we prove Lemma 1 and Theorem 3. Theorems 4 and 5 are proved in Section 4.

## 3. Homogeneous invariant laws.

3.1. *Generalities.* In order to prepare for the proof of Theorem 3, we will derive some results for a general class of cancellative spin systems on $\mathbb{Z}^d$. Throughout this section, $X$ is a cancellative spin system on $\mathbb{Z}^d$ defined by rates $a(A) = a_X(A)$ (see Section 1.2), and $Y$ is its dual cancellative spin system defined by rates $a_Y(A) = a_X(A^T)$. In order for $X$ and $Y$ to be well-defined, we make the summability assumptions:

(3.1)
$$\text{(i)} \quad \sup_i \sum_{A \ni i} a(A)|\{j : (i,j) \in A\}| < \infty,$$
$$\text{(ii)} \quad \sup_i \sum_{A \ni i} a(A)|\{j : (j,i) \in A\}| < \infty.$$

We also assume that our rates are (spatially) homogeneous, in the sense that $a(T_i A) = a(A)$ for all $i \in \mathbb{Z}^d$, and finite $A \subset \mathbb{Z}^d \times \mathbb{Z}^d$, where

(3.2) $\quad T_i A := \{(j+i, k+i) : (j,k) \in A\} \qquad (i, j, k \in \mathbb{Z}^d).$

We will sometimes need the graphical representation of $X$. Independently for each finite $A \in \mathbb{Z}^d \times \mathbb{Z}^d$, let $\pi_A$ be a random, locally finite subset of $\mathbb{R}$, generated by a Poisson processes with intensities $a(A)$. We visualize this by drawing an arrow from $(i,t)$ to $(j,t)$ for each $t \in \pi_A$ and $(j,i) \in A$ (note the order). By definition, a *path* from a subset $C \subset \mathbb{Z}^d \times \mathbb{R}$ to another such subset $D$ is a sequence of points $i_0, \ldots, i_n \in \mathbb{Z}^d$ and times $t_0 \leq t_1 < \cdots < t_n \leq t_{n+1}$ ($n \geq 0$) with

$$(3.3) \qquad \forall 1 \leq m \leq n \; \exists A \text{ s.t. } t_m \in \pi_A \quad \text{and} \quad (i_m, i_{m-1}) \in A,$$

such that $(i_0, t_0) \in C$ and $(i_n, t_{n+1}) \in D$. Thus, a path must walk upward in time and may jump from one site to another along arrows. With these conventions, for any subset $U \subset \mathbb{Z}^d$, setting for $t \geq 0$, $i \in \mathbb{Z}^d$,

$$(3.4) \qquad X_t(i) := 1_{\{\text{the number of paths from } U \times \{0\} \text{ to } (i,t) \text{ is odd}\}}$$

defines a version of the cancellative spin system $X$ defined by the rates $a(X)$, with initial state $X_0 = 1_U$.

PROOF OF LEMMA 1. We only prove part (a), the proof of part (b) being similar. By duality (1.8), for each $y \in \{0,1\}^{\mathbb{Z}^d}$ with $|y| < \infty$,

$$(3.5) \qquad \mathbb{P}[|X_t^{1/2} y| \text{ is odd}] = \mathbb{P}[|X_0^{1/2} Y_t^y| \text{ is odd}] = \tfrac{1}{2}\mathbb{P}[Y_t^y \neq \underline{0}],$$

hence,

$$(3.6) \qquad \mathbb{P}[|X_t^{1/2} y| \text{ is odd}] \underset{t \to \infty}{\longrightarrow} \tfrac{1}{2}\mathbb{P}[Y_t^y \neq \underline{0} \; \forall t \geq 0].$$

Recall that the functions $x \mapsto 1_{\{|xy| \text{ is odd}\}}$ with $|y| < \infty$ are distribution determining. Therefore, the odd upper invariant law of $X$ is uniquely determined by (2.4). If $Y$ survives, then (2.4) shows that $\mathbb{P}[|X_\infty^{1/2} y| \text{ is odd}] > 0$ for some $|y|$ even, hence, $\nu_X^{1/2}$ is not concentrated on $\{\underline{0}, \underline{1}\}$. This shows that (iii)$\Rightarrow$(ii). To show that (ii)$\Rightarrow$(i), it suffices to note that $\mathbb{P}[X_\infty^{1/2} \in \cdot \mid X_\infty^{1/2} \neq \underline{0}, \underline{1}]$ is a coexisting invariant law. To see that (i)$\Rightarrow$(iii), assume that $\mathbb{P}[X_\infty \in \cdot]$ is a coexisting invariant law for $X$. Let $\delta_i \in \{0,1\}^{\mathbb{Z}^d}$ be defined by $\delta_i(j) = 1$ if $i = j$ and $0$ otherwise. Then by the duality (1.8) applied to the process $X$ started in $X_\infty$, we have, for $i \neq j, t \geq 0$,

$$(3.7) \quad \mathbb{P}[Y_t^{\delta_i + \delta_j} \neq \underline{0}] \geq \mathbb{P}[|X_\infty Y_t^{\delta_i + \delta_j}| \text{ is odd}] = \mathbb{P}[X_\infty(i) \neq X_\infty(j)] > 0,$$

and therefore,

$$(3.8) \quad \mathbb{P}[Y_t^{\delta_i + \delta_j} \neq \underline{0} \; \forall t \geq 0] = \lim_{t \to \infty} \mathbb{P}[Y_t^{\delta_i + \delta_j} \neq \underline{0}] \geq \mathbb{P}[X_\infty(i) \neq X_\infty(j)] > 0,$$

which shows that $Y$ survives. $\square$

For later use, we introduce some more notation. The set $\{0,1\}$, equipped with multiplication and addition modulo 2, is a finite field. We can view $\{0,1\}^{\mathbb{Z}^d}$ as a linear space over this field. In this point of view, if we identify a finite set $A \subset \mathbb{Z}^d \times \mathbb{Z}^d$ with the matrix $A$ such that $A(i,j) = 1$ if $(i,j) \in A$ and $A(i,j) = 0$ otherwise, then $Ax$ as defined in (1.4) is the usual action of a matrix on a vector:

$$(3.9) \qquad Ax(i) = \sum_j A(i,j) x(j) \mod(2).$$

In line with these observations, for $x \in \{0,1\}^{\mathbb{Z}^d}$, we let $x^{\mathrm{T}}$ denote the $[\mathrm{mod}(2)]$ linear form on $\{0,1\}^{\mathbb{Z}^d}$ given by

$$(3.10) \qquad x^{\mathrm{T}} y := \sum_i x(i) y(i) \mod(2) = 1_{\{|xy| \text{ is odd}\}},$$

which is well defined whenever $|xy| < \infty$.

3.2. *Uniqueness and convergence.* In this section we continue to work in the general set-up introduced above. We show how a sort of "extinction versus unbounded growth" for the dual process $Y$ can be used to prove convergence to equilibrium for $X$, started in a homogeneous initial law. Our result is similar in spirit to claims by Simonelli [22]. Unfortunately, as already mentioned in Section 2.2, that reference contains an error, which we point out in Section 3.6 below.

To formulate our result, we need to identify "good" configurations where parity can change. To this aim, we select a finite set $\mathcal{B}$ whose elements are finite subsets $B$ of $\mathbb{Z}^d \times \mathbb{Z}^d$ such that $a(B) > 0$, and we define, for $y \in \{0,1\}^{\mathbb{Z}^d}$,

$$(3.11) \quad \|y\|_{\mathcal{B}} := |\{i \in \mathbb{Z}^d : \exists x \in \{0,1\}^{\mathbb{Z}^d}, B \in \mathcal{B} \text{ s.t. } y^{\mathrm{T}}(T_i B) x = 1\}|.$$

Note that $y^{\mathrm{T}}(T_i B) x = 1$ is equivalent to $|\{(j,k) \in T_i B : y(j) = 1 = x(k)\}|$ being odd, or, equivalently, to $|x'y|$ having a different parity from $|xy|$, where $x' = x + (T_i B) x \mod(2)$.

For our models of interest, we may choose $\|\cdot\|_{\mathcal{B}}$ as follows. If $X$ is a neutral Neuhauser–Pacala model affine voter model or rebellious voter model (with arbitrary $\alpha$), then we can find $i,j$ such that

$$(3.12) \qquad a(\{0\} \times \{i,j\}) > 0.$$

Taking for $\mathcal{B}$, the one-point set $\{\{0\} \times \{i,j\}\}$ now leads to $\|y\|_{\mathcal{B}} = |y|$. If, on the other hand, $X$ is the dual of any of these models, then for each $k = 1, \ldots, d$, we can find $i_k$ such that

$$(3.13) \qquad a(\{0, e_k\} \times \{i_k\}) > 0,$$

where $e_1 := (1, 0, \ldots, 0), \ldots, e_d := (0, \ldots, 0, 1)$ are unit vectors in each of the $d$ dimensions of $\mathbb{Z}^d$. Taking for $\mathcal{B}$, the set $\{\{0, e_1\} \times \{i_1\}, \ldots, \{0, e_d\} \times \{i_d\}\}$ now yields $\|y\|_\mathcal{B} = |\nabla y|$ [recall (2.10)].

We say that a probability law $\mu$ on $\{0,1\}^{\mathbb{Z}^d}$ is *locally nonsingular* if, for every finite $\Delta \subset \mathbb{Z}^d$ and every $y \in \{0,1\}^\Delta$,

(3.14) $$\mu(\{x : x(i) = y(i) \ \forall i \in \Delta\}) > 0,$$

that is, every finite configuration has positive probability. We say that $\mu$ is *$X$-nontrivial* if $\mu$ is concentrated on states $x$ such that $\mathbb{P}[X_t^x \in \cdot]$ is locally nonsingular for each $t > 0$.

THEOREM 6 (Homogeneous invariant laws). *Assume that $X$ has a homogeneous $X$-nontrivial invariant law. If the dual process $Y^y$ started in any finite initial state $Y_0^y = y$ satisfies*

(3.15) $\quad \mathbb{P}[\exists t \geq 0 \ s.t. \ \|Y_t^y\|_\mathcal{B} \notin \{1, \ldots, N\}] = 1 \qquad (N \geq 1, \ |y| < \infty),$

*then $\nu_\mathrm{X}^{1/2}$ is $X$-nontrivial, and the unique homogeneous $X$-nontrivial invariant law of $X$. If, moreover, one has*

(3.16) $\quad \lim_{t \to \infty} \mathbb{P}[\|Y_t^y\|_\mathcal{B} \notin \{1, \ldots, N\}] = 1 \qquad (N \geq 1, |y| < \infty),$

*then the process $X$ started in any homogeneous $X$-nontrivial initial law satisfies*

(3.17) $$\mathbb{P}[X_t \in \cdot] \underset{t \to \infty}{\Longrightarrow} \nu_\mathrm{X}^{1/2}.$$

The proof depends on two lemmas.

LEMMA 7 (Parity uncertainty). *For each $\varepsilon > 0$ and $t > 0$, there exists an $N \geq 1$ such that if $X$ and $Y$ are started in deterministic initial states $X_0$ and $Y_0$ satisfying*

(3.18) $$|\{i \in \mathbb{Z}^d : \exists B \in \mathcal{B} \ s.t. \ Y_0^\mathrm{T}(T_i B) X_0 = 1\}| \geq N,$$

*then*

(3.19) $$|\mathbb{P}[|X_t Y_0| \ is \ odd] - \tfrac{1}{2}| \leq \varepsilon.$$

PROOF. Let $\mathcal{A}$ denote the set of all finite subsets of $\mathbb{Z}^d \times \mathbb{Z}^d$. For any $A \in \mathcal{A}$, let us put

(3.20) $$\begin{aligned} R_-(A) &:= \{j : A(i,j) = 1 \text{ for some } i \in \mathbb{Z}^d\}, \\ R_+(A) &:= \{i : A(i,j) = 1 \text{ for some } j \in \mathbb{Z}^d\}, \end{aligned}$$

and let us call $R(A) := R_-(A) \cup R_+(A)$ the *range* of $A$. Let us say that $A, \tilde{A} \in \mathcal{A}$ are *disjoint* if $R(A) \cap R(\tilde{A}) = \varnothing$. If $X_0, Y_0$ satisfy (3.18), then it is not hard to see that we can successively choose $i_1, \ldots, i_n \in \mathbb{Z}^d$ and $B_1, \ldots, B_n \in \mathcal{B}$ such that $T_{i_1} B_1, \ldots, T_{i_n} B_n$ are disjoint, $Y_0^{\mathrm{T}}(T_{i_m} B_m) X_0 = 1$ for all $1 \leq m \leq n$, and

(3.21) $$n \geq N/(|\mathcal{B}| K^2),$$

where $K := \max\{|R(B)| : B \in \mathcal{B}\}$. Indeed, for each $B, B' \in \mathcal{B}$ and $i \in \mathbb{Z}^d$, there are at most $K^2$ points $j \in \mathbb{Z}^d$ such that $R(T_j B') \cap R(T_i B) \neq \varnothing$, so once we pick a point from the set in (3.18), there are at most $|\mathcal{B}| K^2$ points we cannot pick anymore. We now use the graphical representation (see Section 3.1). Let

(3.22) $$\begin{aligned} M := \{m : 1 \leq m \leq n, |\pi_{T_{i_m} B_m} \cap (0,t)| \in \{0,1\} \text{ and} \\ |\pi_A \cap (0,t)| = 0 \ \forall A \in \mathcal{A}, A \neq T_{i_m} B_m \\ \text{s.t. } R(A) \cap R(T_{i_m} B_m) \neq \varnothing\}. \end{aligned}$$

Thus, for each $m \in M$, the Poisson process associated with $T_{i_m} B_m$ becomes active zero or one time during the time interval $(0, t)$, and no other Poisson process creates arrows in the range $R(T_{i_m} B_m)$ during this time interval. We claim that if $N$ is sufficiently large, then the set $M$ is large with high probability.

Indeed, since $a(T_{i_m} B_m) = a(B_m)$ is bounded from above, the probability that $|\pi_{T_{i_m} B_m} \cap (0,t)| \in \{0,1\}$ is uniformly bounded from below. By summability [condition (3.1)], the probability that no other Poisson process creates arrows in the range $R(T_{i_m} B_m)$ during $(0, t)$ is also uniformly bounded from below. These events are not independent for different $m$, but they are positively correlated, so we get a lower bound assuming independence, which proves our claim.

Let $\mathcal{F}$ be the $\sigma$-field generated by the random set $M$ and by all Poisson processes $\pi_A \cap (0, t)$ with $A \in \mathcal{A} \setminus \{T_{i_m} B_m : m \in M\}$. Thus, $\mathcal{F}$ corresponds to knowing the random set $M$ and all Poisson processes on $(0, t)$, except those associated with the $T_{i_m} B_m$ with $m \in M$. Set $\theta_m := a(B_m) t$. Note that the $\theta_m$ are uniformly bounded from above and below by the fact that $\mathcal{B}$ is finite. We claim that if we condition on $\mathcal{F}$, then under the conditioned law, the random variables $|\pi_{T_{i_m} B_m} \cap (0,t)|$ are independent $\{0, 1\}$-valued random variables with

(3.23) $$\mathbb{P}[|\pi_{T_{i_m} B_m} \cap (0,t)| = 1 \mid \mathcal{F}] = \frac{\theta_m e^{-\theta_m}}{e^{-\theta_m} + \theta_m e^{-\theta_m}} =: \phi_m \qquad (m \in M).$$

Indeed, first condition on all $\pi_A \cap (0, t)$ with $A \neq T_{i_1} B_1, \ldots, T_{i_n} B_n$. Under this conditional law, the $\pi_{T_{i_1} B_1} \cap (0, t), \ldots, \pi_{T_{i_n} B_n} \cap (0, t)$ are independent

Poisson processes with intensities $a(B_1),\ldots,a(B_n)$. Let

$$(3.24) \quad M' := \{m : 1 \leq m \leq n, |\pi_A \cap (0,t)| = 0 \ \forall A \in \mathcal{A}\setminus\{T_{i_1}B_1,\ldots,T_{i_n}B_n\}$$
$$\text{s.t. } R(A) \cap R(T_{i_m}B_m) \neq \varnothing\}.$$

Under the conditional law we are considering, $M'$ is a deterministic set, and the $\pi_{T_{i_m}B_m} \cap (0,t)$ with $m \in M'$ are independent Poisson processes. Hence, since

$$(3.25) \quad M = \{m \in M' : |\pi_{T_{i_m}B_m} \cap (0,t)| \in \{0,1\}\},$$

if we condition also on the $\pi_{T_{i_m}B_m} \cap (0,t)$ with $m \in M'\setminus M$, then under this new conditional law, the $\pi_{T_{i_m}B_m} \cap (0,t)$ with $m \in M$ are independent Poisson processes conditioned to produce zero or one point. This explains (3.23).

Using the graphical representation, we now write $X_0 = 1_U$, $Y_0 = 1_V$, and

$$(3.26) \quad \mathbb{P}[|X_tY_0| \text{ is even}] - \mathbb{P}[|X_tY_0| \text{ is odd}] = \mathbb{E}[(-1)^{P+P'}],$$

where $P$ and $P'$ are the number of paths from $U \times \{0\}$ to $V \times \{t\}$ that do and do not use, respectively, arrows created by the Poisson processes $\pi_{T_{i_m}B_m} \cap (0,t)$ with $m \in M$. Note that due to the definition of $M$ the paths counted by $P$ use exactly one arrow on $(0,t)$ created by a Poisson processes $\pi_{T_{i_m}B_m} \cap (0,t)$ with $m \in M$. Since $Y_0^T(T_{i_m}B_m)X_0 = 1$ for $m \in M$, this implies that $P = \sum_{m \in M} |\pi_{T_{i_m}B_m} \cap (0,t)| \bmod(2)$. It therefore follows that

$$(3.27) \quad \begin{aligned} \mathbb{E}[(-1)^{P+P'} \mid \mathcal{F}] &= (-1)^{P'}\mathbb{E}[(-1)^P \mid \mathcal{F}] \\ &= (-1)^{P'} \prod_{m \in M} \mathbb{E}[(-1)^{|\pi_{T_{i_m}B_m} \cap (0,t)|} \mid \mathcal{F}] \\ &= (-1)^{P'} \prod_{m \in M} ((1-\phi_m) - \phi_m), \end{aligned}$$

where the $\phi_m$ are defined in (3.23). Integrating over the $\sigma$-field $\mathcal{F}$, it follows that, under the unconditional law,

$$(3.28) \quad |\mathbb{P}[|X_tY_0| \text{ is even}] - \mathbb{P}[|X_tY_0| \text{ is odd}]| \leq \mathbb{E}\left[\prod_{m \in M} |1 - 2\phi_m|\right].$$

Since the $\phi_m$ are bounded away from zero and one and since $|M|$ is with high probability large if $N$ is large, it follows that for each $\varepsilon > 0$ we can choose $N$ large enough such that (3.19) holds. $\square$

LEMMA 8 (Many good configurations). *Let $X$ be started in a homogeneous $X$-nontrivial initial law and let $t > 0$. Assume that $y_n \in \{0,1\}^{\mathbb{Z}^d}$ satisfy*

$$(3.29) \quad \lim_{n \to \infty} \|y_n\|_{\mathcal{B}} = \infty.$$

*Then*

$$(3.30) \qquad |\{i \in \mathbb{Z}^d : \exists B \in \mathcal{B} \text{ s.t. } y_n^{\mathrm{T}}(T_i B) X_t = 1\}| \xrightarrow[n \to \infty]{\mathrm{P}} \infty,$$

*where $\xrightarrow{\mathrm{P}}$ denotes convergence in probability.*

PROOF. This proof very closely follows ideas from [15], Theorem (9.2). Set

$$(3.31) \quad \begin{aligned} C_n &:= \{i \in \mathbb{Z}^d : \exists x \in \{0,1\}^{\mathbb{Z}^d}, \ B \in \mathcal{B} \text{ s.t. } y_n^{\mathrm{T}}(T_i B) x = 1\}, \\ C_n' &:= \{i \in \mathbb{Z}^d : \exists B \in \mathcal{B} \text{ s.t. } y_n^{\mathrm{T}}(T_i B) X_t = 1\}. \end{aligned}$$

By (3.29), $|C_n| \to \infty$. We need to show that the random subsets $C_n' \subset C_n$ satisfy

$$(3.32) \qquad \mathbb{P}[|C_n'| \geq N] \xrightarrow[n \to \infty]{} 1$$

for all $N \geq 1$. By dividing $C_n$ into $N$ disjoint sets, each with size tending to infinity, we can reduce this to showing (3.32) for $N = 1$. For each $i \in C_n$, choose $B_{i,n} \in \mathcal{B}$ and $x_{i,n} \in \{0,1\}^{\mathbb{Z}^d}$ such that

$$(3.33) \qquad y_n^{\mathrm{T}}(T_i B_{i,n}) x_{i,n} = 1.$$

Since $\mathcal{B}$ is finite, by going to a subsequence if necessary, we can assume without loss of generality that for some $B \in \mathcal{B}$ and $z \in \{0,1\}^{R_+(B)}$, where $R_+(B)$ is as in (3.20), $|\tilde{C}_n| \to \infty$, where

$$(3.34) \qquad \tilde{C}_n := \{i \in C_n : B_{i,n} = B, (x_{i,n}(i+j))_{j \in R_+(B)} = z\}.$$

It now suffices to prove that $\mathbb{P}[|\tilde{C}_n'| \geq 1] \to 1$, where

$$(3.35) \qquad \tilde{C}_n' := \{i \in C_n : (X_t(i+j))_{j \in R_+(B)} = z\}.$$

Equivalently, we need to show that

$$(3.36) \qquad \mathbb{E}\left[\prod_{i \in \tilde{C}_n} 1_{\{(X_t(i+j))_{j \in R_+(B)} \neq z\}}\right] \xrightarrow[n \to \infty]{} 0.$$

Fix $\varepsilon > 0$ and $k \geq 1$. For each $L \geq 1$, we can find, for $n$ sufficiently large, subsets $\tilde{C}_n^{L,k} \subset \tilde{C}_n$ such that $|i-j| \geq L$ for all $i, j \in \tilde{C}_n^{L,k}$, $i \neq j$ and $|\tilde{C}_n^{L,k}| = k$. We claim that there exists an $L \geq 1$ such that the process $X^x$ started in any deterministic initial state $X_0^x = x \in \{0,1\}^{\mathbb{Z}^d}$ satisfies

$$(3.37) \quad \begin{aligned} E&\left[\prod_{i \in \tilde{C}_n^{L,k}} 1_{\{(X_t^x(i+j))_{j \in R_+(B)} \neq z\}}\right] \\ &\leq \prod_{i \in \tilde{C}_n^{L,k}} \mathbb{E}[1_{\{(X_t^x(i+j))_{j \in R_+(B)} \neq z\}}] + \varepsilon k. \end{aligned}$$

One way to see this is to check that the conditions of [17], Theorem I.4.6 are fulfilled. Alternatively, one can use the graphical representation. It is not hard to see that, for $L$ sufficiently large, the probability that there exist two paths between time zero and time $t$, one ending at $R_+(T_iB)$ and the other at $R_+(T_jB)$ for some $i,j \in \tilde{C}_n^{L,k}$, $i \neq j$, and both starting at the same site, is bounded by $\varepsilon k$. This implies (3.37).

By (3.37), Hölder's inequality and the fact that $\mathbb{P}[X_0 \in \cdot]$ is homogeneous, it follows that

$$\limsup_{n\to\infty} \mathbb{E}\left[\prod_{i\in\tilde{C}_n} 1_{\{(X_t(i+j))_{j\in R_+(B)}\neq z\}}\right]$$

$$\leq \limsup_{n\to\infty} \int \mathbb{P}[X_0 \in dx]\mathbb{E}\left[\prod_{i\in\tilde{C}_n^{L,k}} 1_{\{(X_t^x(i+j))_{j\in R_+(B)}\neq z\}}\right]$$

(3.38) $$\leq \limsup_{n\to\infty} \int \mathbb{P}[X_0 \in dx] \prod_{i\in\tilde{C}_n^{L,k}} \mathbb{E}[1_{\{(X_t^x(i+j))_{j\in R_+(B)}\neq z\}}] + \varepsilon k$$

$$\leq \limsup_{n\to\infty} \prod_{i\in\tilde{C}_n^{L,k}} \left(\int \mathbb{P}[X_0 \in dx]\mathbb{P}[(X_t^x(i+j))_{j\in R_+(B)} \neq z]^k\right)^{1/k}$$

$$+ \varepsilon k$$

$$= \int \mathbb{P}[X_0 \in dx]\mathbb{P}[(X_t^x(j))_{j\in R_+(B)} \neq z]^k + \varepsilon k.$$

Letting first $\varepsilon \to 0$ and then $k \to \infty$, using nontriviality, we arrive at (3.36).  □

Lemmas 7 and 8 combine to give the following corollary.

COROLLARY 9 (Parity indeterminacy). *Let $X$ be started in a homogeneous $X$-nontrivial initial law and let $t > 0$. Assume that $y_n \in \{0,1\}^{\mathbb{Z}^d}$ satisfy*

(3.39) $$\lim_{n\to\infty} \|y_n\|_{\mathcal{B}} = \infty.$$

*Then*

(3.40) $$\lim_{n\to\infty} \mathbb{P}[|X_ty_n| \text{ is odd}] = \tfrac{1}{2}.$$

PROOF OF THEOREM 6. Imagine that $\nu$ is a homogeneous $X$-nontrivial invariant law of $X$. Put

(3.41) $$f(y) := \int \nu(dx) 1_{\{|xy| \text{ is odd}\}} \qquad (y \in \{0,1\}^{\mathbb{Z}^d},\ |y| < \infty).$$

We claim that $(f(Y_t^y))_{t\geq 0}$ is a martingale. Indeed, if $X$ is a stationary process with law $\mathbb{P}[X_t \in \cdot] = \nu$ $(t \in \mathbb{R})$, independent of $Y^y$, then, by duality (1.8) for $y \in \{0,1\}^{\mathbb{Z}^d}$ with $|y| < \infty$,

$$(3.42) \quad \mathbb{E}[f(Y_t^y)] = \mathbb{P}[|X_0 Y_t^y| \text{ is odd}] = \mathbb{P}[|X_t Y_0^y| \text{ is odd}] = f(y).$$

Using, moreover, the Markov property of $Y$, this shows that $(f(Y_t^y))_{t\geq 0}$ is a bounded martingale. Set

$$(3.43) \qquad \tau_N := \inf\{t \geq 0 : \|Y_t^y\|_{\mathcal{B}} \notin \{1, \ldots, N\}\},$$

which is a.s. finite for all $N \geq 1$ by our assumption (3.15). Hence, by optional stopping,

$$(3.44) \quad f(y) = \mathbb{E}[f(Y_{\tau_N}^y)] = \mathbb{P}[|X_t Y_{\tau_N}^y| \text{ is odd} | \|Y_{\tau_N}^y\|_{\mathcal{B}} > N] \mathbb{P}[\|Y_{\tau_N}^y\|_{\mathcal{B}} > N].$$

Letting $N \to \infty$, using Corollary 9, we find that for $y \in \{0,1\}^{\mathbb{Z}^d}$ with $|y| < \infty$,

$$(3.45) \qquad \int \nu(dx) 1_{\{|xy| \text{ is odd}\}} = f(y) = \tfrac{1}{2} \mathbb{P}[Y_t^y \neq \underline{0} \ \forall t \geq 0].$$

By (2.4), this implies that $\nu = \nu_X^{1/2}$.

If, moreover, (3.16) holds and $X$ is started in any homogeneous $X$-nontrivial initial law, then, by duality and Corollary 9,

$$(3.46) \quad \mathbb{P}[|X_t y| \text{ is odd}] = \mathbb{P}[|X_s Y_{t-s}^y| \text{ is odd}] \xrightarrow[t\to\infty]{} \tfrac{1}{2} \mathbb{P}[Y_t^y \neq \underline{0} \ \forall t \geq 0].$$

Since this holds for any finite $y$, it follows that $\mathbb{P}[X_t \in \cdot] \Rightarrow \nu_X^{1/2}$ as $t \to \infty$. $\square$

3.3. *Extinction versus unbounded growth.* In order for Theorem 6 to be applicable, we need to check that the dual $Y$ of a cancellative spin system $X$ satisfies a version of "extinction versus unbounded growth," that is, we must show that $\|Y_t^y\|_{\mathcal{B}}$ is either zero or large at random times $t$ or at large fixed $t$. Therefore, in this section, we derive sufficient conditions for a cancellative spin system to show this kind of behavior.

We start with some simple observations. Let $X$ be a nonexplosive continuous-time Markov process with countable state space $S$ (e.g., a cancellative spin system restricted to the space of finite states). Let $X^x$ denote the process $X$ started in $X_0^x = x$ and let $D \subset S$.

LEMMA 10 (Markov process leaving sets). (a) *If* $\inf_{x \in D} \mathbb{P}[\exists t \geq 0 \ s.t. \ X_t^x \notin D] > 0$, *then* $\mathbb{P}[\exists t \geq 0 \ s.t. \ X_t^x \notin D] = 1$ *for all* $x \in S$.

(b) *If* $\inf_{x \in D} \mathbb{P}[\exists t \geq 0 \ s.t. \ X_u^x \notin D \ \forall u \geq t] > 0$, *then* $\mathbb{P}[\exists t \geq 0 \ s.t. \ X_u^x \notin D \ \forall u \geq t] = 1$ *for all* $x \in S$.

PROOF. To prove part (b), set $\varepsilon := \inf_{x \in D} \mathbb{P}[\exists t \geq 0 \text{ s.t. } X_u^x \notin D \; \forall u \geq t]$. Let $(\mathcal{F}_t)_{t \geq 0}$ be the filtration generated by $X^x$. By the Markov property and martingale convergence,

$$
\begin{aligned}
\varepsilon 1_D(X_s) &\leq \mathbb{P}[\exists t \geq 0 \text{ s.t. } X_u \notin D \; \forall u \geq t \mid X_s] \\
&= \mathbb{P}[\exists t \geq 0 \text{ s.t. } X_u \notin D \; \forall u \geq t \mid \mathcal{F}_s] \\
&\xrightarrow[s \to \infty]{} 1_{\{\exists t \geq 0 \text{ s.t. } X_u \notin D \; \forall u \geq t\}} \quad \text{a.s.}
\end{aligned}
$$
(3.47)

This shows that $\lim_{s \to \infty} 1_D(X_s) = 0$ a.s. on the complement of the event $\{\exists t \geq 0 \text{ s.t. } X_u \notin D \; \forall u \geq t\}$, which implies our claim. Part (a) follows from part (b), applied to the process stopped at $\tau := \inf\{t \geq 0 : X_t^x \notin D\}$. □

As a simple consequence of Lemma 10, we obtain the following corollary. Recall the definition of $|\nabla X|$ in (2.10).

COROLLARY 11 (Extinction versus unbounded growth). *Let $\alpha < 1$ and let $X$ be either the neutral Neuhauser–Pacala model, affine voter model or rebellious voter model, and let $Y$ be its dual parity preserving branching process. Then:*

(a) *One has*

$$
\begin{aligned}
\mathbb{P}[\exists t \geq 0 \text{ s.t. } |\nabla X_t^x| \notin \{1, \ldots, N\}] &= 1 \quad (N \geq 1, \; |x| < \infty), \\
\mathbb{P}[\exists t \geq 0 \text{ s.t. } |Y_t^y| \notin \{1, \ldots, N\}] &= 1 \quad (N \geq 1, \; |y| < \infty).
\end{aligned}
$$
(3.48)

(b) *If, moreover, $d \geq 2$ and $\alpha > 0$, then*

$$
(3.49) \quad \mathbb{P}\left[\lim_{t \to \infty} |\nabla X_t^x| = \infty \text{ or } \exists t \geq 0 \text{ s.t. } X_t^x = \underline{0}\right] = 1 \quad (|x| < \infty).
$$

PROOF. We claim that, for all $N \geq 1$ and $t > 0$,

$$
\begin{aligned}
&\text{(i)} \quad \inf_{|\nabla x| \leq N} \mathbb{P}[|\nabla X_t^x| \notin \{1, \ldots, N\}] > 0, \\
&\text{(ii)} \quad \inf_{|y| \leq N} \mathbb{P}[|Y_t^y| \notin \{1, \ldots, N\}] > 0,
\end{aligned}
$$
(3.50)

and, if $\alpha > 0$,

$$
(3.51) \quad \inf_{0 < |x| \leq N} \mathbb{P}[X_t^x = \underline{0}] > 0.
$$

Indeed, (3.50)(ii) follows from the fact that, by our assumption that $\alpha < 1$, a particle lying sufficiently on the "outside" of $y$ may produce $N$ more particles before anything happens to the other particles of $y$. Formula (3.50)(i) follows from similar considerations. If $\alpha > 0$, then by the voter model dynamics, a collection of at most $N$ ones has a uniformly positive probability to die out in time $t$, which proves (3.51).

By Lemma 10(a), (3.50) implies (3.48). By Lemma 10(b) applied to sets of the form $D := \{x \in \{0,1\}^{\mathbb{Z}^d} : 0 < |x| \leq N\}$, (3.51) implies that $\mathbb{P}[A_N^x] = 1$ for any $N \geq 1, |x| < \infty$, where

(3.52) $\qquad A_N^x := \{\exists t \geq 0 \text{ s.t. } |X_u^x| \geq N \text{ or } X_t^x = \underline{0} \ \forall u \geq t\}.$

Since $A_N^x \downarrow A^x$ with

(3.53) $\qquad A^x := \left\{ \lim_{t \to \infty} |X_t^x| = \infty \text{ or } \exists t \geq 0 \text{ s.t. } X_t^x = \underline{0} \right\},$

we obtain $\mathbb{P}[A^x] = 1$. In dimensions $d \geq 2$ this implies (3.49). $\square$

In order to prove Theorem 3(a), we need to prove extinction versus unbounded growth for the parity preserving branching process $Y$. In this case, Lemma 10(b) is of no use, since the analogue of (3.51) for $Y$ does not hold because of parity preservation. In the following Section 3.4 we give a proof assuming that $Y$ is not stable and $\alpha > 0$. An alternative approach, that works only for very small $\alpha$, but includes $\alpha = 0$, is to use comparison with oriented percolation. For the rebellious voter model, we will use this approach in Section 4.

3.4. *Instability.* The main result of this section is the next theorem, which will be used in the proof of Theorem 3(a). We will also apply this result in [23], which is why we formulate it in some generality here.

THEOREM 12 (Extinction versus unbounded growth for parity preserving branching). *Let $Y$ be a spatially homogeneous parity preserving cancellative spin system on $\mathbb{Z}^d$ defined by rates $a(A)$ satisfying (3.1). Assume that for $L \geq 1, n \geq 0, t > 0$,*

(3.54) $\quad \inf\{\mathbb{P}[|Y_t^y| = n] : |y| = n+2, y(i) = 1 = y(j)$

$\qquad \qquad \qquad \qquad \textit{for some } i \neq j, |i-j| \leq L\} > 0.$

*Assume that $Y$ is not stable. Then*

(3.55) $\qquad \lim_{t \to \infty} \mathbb{P}[0 < |Y_t^y| < N] = 0 \qquad (N \geq 1, |y| < \infty).$

As a first step, we prove that the convergence in (3.55) holds in Césaro mean.

PROPOSITION 13 (Extinction versus unbounded growth). *Under the same assumptions as in Theorem 12,*

(3.56) $\qquad \lim_{T \to \infty} \frac{1}{T} \int_0^T dt\, \mathbb{P}[0 < |Y_t^y| < N] = 0 \qquad (N \geq 0, |y| < \infty).$

PROOF. The idea of the proof is as follows. By our assumption that $Y$ is not stable, one particle alone will soon produce at least three particles. In fact, if the process does not die out, then most of the time it will contain at least three particles. These particles cannot stay close together, for else they would annihilate each other. But single particles far from each other will soon each again produce at least three particles, and therefore, the number of particles must keep growing.

To make this idea precise, we use induction. We write $|y| = n \mod(2)$ to indicate that $|y| < \infty$, and $|y|$ is even or odd depending on whether $n$ is even or odd. For $n \geq 0$, consider the following statements:

$$(\mathrm{I}_n) \quad \lim_{T\to\infty} \frac{1}{T} \int_0^T dt\, \mathbb{P}[0 < |Y_t^y| \leq n] = 0 \qquad \text{for all } |y| = n \mod(2).$$

$$(\mathrm{II}_n) \quad \lim_{T\to\infty} \frac{1}{T} \int_0^T dt\, \mathbb{P}[\exists s \in [0, S] \text{ s.t. } 0 < |Y_{t+s}^y| \leq n] = 0$$
$$\text{for all } |y| = n \mod(2), S > 0.$$

$$(\mathrm{III}_n) \quad \lim_{T\to\infty} \frac{1}{T} \int_0^T dt\, \mathbb{P}[\exists s \in [0, S] \text{ s.t. } |Y_{t+s}^y| = n+2 \text{ and}$$
$$Y_{t+s}^y(i) = 1 = Y_{t+s}^y(j) \qquad \text{for some } i \neq j,\ |i-j| \leq L] = 0$$
$$\text{for all } |y| = n \mod(2), S > 0, L \geq 1.$$

$$(\mathrm{IV}_n) \quad \lim_{S\to\infty} \limsup_{T\to\infty} \mathbb{E}\left[\frac{1}{T}\int_0^T dt\, \frac{1}{S}\int_0^S ds\, 1_{\{\tau_n(t) \leq S, 0 < |Y_{t+\tau_n(t)+s}^y| \leq n\}}\right] = 0$$
$$\text{for all } |y| = n \mod(2), \text{ where } \tau_n(t) := \inf\{u \geq 0 : 0 < |Y_{t+u}^y| \leq n\}.$$

$$(\mathrm{V}_n) \quad \lim_{S\to\infty} \limsup_{T\to\infty} \frac{1}{T}\int_0^T dt\, \frac{1}{S}\int_0^S ds\, \mathbb{P}[0 < |Y_{t+s}^y| \leq n] = 0$$
$$\text{for all } |y| = n \mod(2).$$

We will prove that $\mathrm{I}_0$ and $\mathrm{I}_1$ hold, and $\mathrm{I}_n$ implies $\mathrm{I}_{n+2}$. Observe that if $\sigma_T$ is uniformly distributed on $[0, T]$ and independent of $Y^y$, then

$$(3.57) \qquad \mathbb{P}[0 < |Y_{\sigma_T}^y| \leq n] = \mathbb{E}\left[\frac{1}{T}\int_0^T dt\, 1_{\{0 < |Y_t^y| \leq n\}}\right].$$

In the proofs below we will freely change between these and similar formulas.

$\mathrm{I}_0$ *and* $\mathrm{I}_1$ *hold.* $\mathrm{I}_0$ is trivial. Since, by assumption, $Y$ viewed from its lower left corner is not positively recurrent, the probability that $Y_t$ consists of a single particle tends to zero as $t \to \infty$, which proves $\mathrm{I}_1$.

$\mathrm{I}_n$ *implies* $\mathrm{II}_n$. This follows from the observation that

$$(3.58) \qquad \inf\{\mathbb{P}[0 < |Y_t^y| \leq n\ \forall 0 \leq t \leq 1] : 0 < |y| \leq n\} =: p > 0.$$

Indeed, conditional on $0 < |Y^y_{t+s}| \leq n$ for some $s \in [0, S]$, with probability at least $p$, the process has between 1 and $n$ particles during a time interval of length one somewhere between time $t$ and $t + S + 1$. Therefore, if $\sigma_T$ and $\sigma_{S+1}$ are uniformly distributed on $[0, T]$ and $[0, S + 1]$, respectively, independent of each other and of $Y^y$, then

$$\limsup_{T \to \infty} \mathbb{P}[\exists s \in [0, S] \text{ s.t. } 0 < |Y^y_{\sigma_T + s}| \leq n]$$

(3.59)
$$\leq \frac{S+1}{p} \limsup_{T \to \infty} \mathbb{P}[0 < |Y^y_{\sigma_T + \sigma_{S+1}}| \leq n]$$

$$= \frac{S+1}{p} \limsup_{T \to \infty} \mathbb{P}[0 < |Y^y_{\sigma_T}| \leq n] = 0,$$

where in the last two steps we have used that the total variation distance between $\mathbb{P}[\sigma_T \in \cdot]$ and $\mathbb{P}[\sigma_T + \sigma_{S+1} \in \cdot]$ tends to zero as $T \to \infty$, and our assumption $\mathrm{I}_n$, respectively.

$\mathrm{II}_n$ *implies* $\mathrm{III}_n$. For $n \geq 1$, using (3.54), it is not hard to see that

(3.60)
$$\inf\{\mathbb{P}[0 < |Y^y_t| \leq n \ \forall 1 \leq t \leq 2] : |y| = n + 2,$$
$$y(i) = 1 = y(j) \text{ for some } i \neq j, |i - j| \leq L\} > 0.$$

From this the implication follows much in the spirit of the previous implication. This argument does not work for $n = 0$, so we will prove that $\mathrm{III}_0$ holds by different means. We observe that for $N = 0$ (3.54) implies that for $L \geq 1$ we obtain for all $t > 0$,

(3.61) $\quad \inf\{\mathbb{P}[Y^y_t = \underline{0}] : y = \delta_i + \delta_j \text{ for some } i \neq j, |i - j| \leq L\} > 0.$

By Lemma 10(b), this implies that, for $L \geq 1$,

(3.62) $\quad \mathbb{P}\left[\lim_{t \to \infty} 1_{\{Y^y_t = \delta_i + \delta_j \text{ for some } i \neq j, |i-j| \leq L\}} = 0\right] = 1,$

which in turn implies $\mathrm{III}_0$.

$\mathrm{I}_1$, $\mathrm{II}_n$ *and* $\mathrm{III}_n$ *imply* $\mathrm{IV}_{n+2}$. Let $\sigma_S$ be uniformly distributed on $[0, S]$, independent of $Y^y$. By $\mathrm{I}_1$ and parity preservation, for each $\varepsilon > 0$, we can choose $S_0 > 0$ such that, for all $S \geq S_0$,

(3.63) $\quad \mathbb{P}[|Y^{\delta_0}_{\sigma_S}| \geq 3] \geq 1 - \varepsilon.$

For each such $S$, we can choose $L \geq 1$ such that

(3.64) $\quad \mathbb{P}\left[Y^{\delta_0}_s(i) = 0 \ \forall 0 \leq s \leq S, |i| \geq \frac{L}{2}\right] \geq 1 - \varepsilon.$

Therefore, if we start the process in a state $y$ such that $|y| = n + 2$ and $|i - j| > L$ for all $i \neq j$ with $y(i) = 1 = y(j)$, then with probability at least

26 A. STURM AND J. SWART

$(1-2\varepsilon)^{n+2}$, all the $n+2$ particles have produced at least 3 particles at time $\sigma_S$, without being influenced by each other. Thus,

$$(3.65) \qquad \mathbb{P}[|Y^y_{\sigma_S}| \geq 3(n+2)] \geq (1-2\varepsilon)^{n+2} =: 1-\varepsilon',$$

for each such $y$. Let $\sigma_T$ be uniformly distributed on $[0,T]$ and independent of $\sigma_S$ and $Y^y$. If $T$ is large, then $\mathrm{II}_n$ and $\mathrm{III}_n$ tell us that the probability that $\tau_{n+2}(\sigma_T) \leq S$ while $Y^y_{\sigma_T + \tau_{n+2}(\sigma_T)}$ does not consist of $n+2$ particles, situated at distance at least $L$ from each other, is small. Therefore, by what we have just proved,

$$(3.66) \qquad \limsup_{T \to \infty} \mathbb{P}[\tau_{n+2}(\sigma_T) \leq S, |Y^y_{\sigma_T + \tau_{n+2}(\sigma_T) + \sigma_S}| \leq n+2] \leq \varepsilon'.$$

Since $\varepsilon'$ can be made arbitrarily small, this proves $\mathrm{IV}_{n+2}$.

$\mathrm{IV}_n$ *implies* $\mathrm{V}_n$. One has

$$\begin{aligned}
\frac{1}{S}\int_0^S ds\, \mathbb{P}[0<|Y^y_{t+s}|\leq n] &= \mathbb{E}\Big[\frac{1}{S}\int_{\tau_n(t)}^S ds\, 1_{\{0<|Y^y_{t+s}|\leq n\}}\Big] \\
(3.67) \qquad &\leq \mathbb{E}\Big[1_{\{\tau_n(t)\leq S\}}\frac{1}{S}\int_{\tau_n(t)}^{\tau_n(t)+S} ds\, 1_{\{0<|Y^y_{t+s}|\leq n\}}\Big] \\
&= \mathbb{E}\Big[\frac{1}{S}\int_0^S ds\, 1_{\{\tau_n(t)\leq S,\ 0<|Y^y_{t+\tau_n(t)+s}|\leq n\}}\Big].
\end{aligned}$$

Integrating from 0 to $T$, dividing by $T$, and taking the limsup as $T \to \infty$ and then the limit $S \to \infty$, the claim follows.

$\mathrm{V}_n$ *implies* $\mathrm{I}_n$. Let $\sigma_T$ and $\sigma_S$ be uniformly distributed on $[0,T]$ and $[0,S]$, respectively, independent of each other and of $Y^y$. By $\mathrm{V}_n$, we can choose $S(T)$ such that $\lim_{T \to \infty} S(T)/T = 0$ and

$$(3.68) \qquad \lim_{T \to \infty} \mathbb{P}[0 < |Y^y_{\sigma_T + \sigma_{(S(T))}}| \leq n] = 0.$$

Since $S(T) \ll T$, the total variation distance between $\mathbb{P}[\sigma_T + \sigma_{S(T)} \in \cdot]$ and $\mathbb{P}[\sigma_T \in \cdot]$ tends to zero as $T \to \infty$, so $\mathrm{I}_n$ follows. $\square$

To get from Proposition 13 to Theorem 12, we need the following lemma, which depends on another lemma.

LEMMA 14 (Aperiodicity). *Under the same assumptions as in Theorem 12,*

$$(3.69) \qquad \lim_{t \to \infty} \sup_{s \in [0,S]} |\mathbb{P}[0<|Y^y_t|<N] - \mathbb{P}[0<|Y^y_{t+s}|<N]| = 0$$

$(N \geq 1, S > 0)$.

PROOF. We need to prove something like aperiodicity for an interacting particle system. This is not an easy problem in general. The only general result that we are aware of is restricted to one-dimensional systems and due to Mountford [20]. In our present setting, however, we can use the fact that each time when there are less than $N$ particles, the next jump of our system happens after an exponential time with mean bounded from below, which causes enough uncertainty in the time variable to prove (3.69). We got this idea from [14], Lemma 2.4. Our way of implementing this idea is quite different from that reference, though.

To make this idea rigorous, we proceed as follows. Let $Z^y$ be the embedded Markov chain associated with $Y^y$, which is defined as follows. For each $y \in \{0,1\}^{\mathbb{Z}^d}$ such that $|y| < \infty$, let

$$r(y) := \sum_A a(A) 1_{\{Ay \neq \underline{0} \mod(2)\}}. \tag{3.70}$$

Note that $r(y)$ is the total rate of jumps of $Y$ from the state $y$ to any other state. It follows from (3.54) that $r(y) > 0$ for all $y \neq \underline{0}$. Now $Z^y = (Z_n^y)_{0 \leq n \leq n_\infty}$ is the Markov chain with state space $\{y \in \{0,1\}^{\mathbb{Z}^d} : |y| < \infty\}$, started in $Z_0^y = y$, which jumps from the state $y$ to the state $y + Ay \mod(2)$ with probability $a(A)/r(y)$. The process $Z^y$ is defined up to the random time $n_\infty := \inf\{n \geq 0 : y = \underline{0}\}$, which may be infinite.

We may construct $Y^y$ from $Z^y$ as follows. Let $(\sigma_n)_{n \geq 0}$ be i.i.d. exponentially distributed random variables with mean one, independent of $Z^y$. Set

$$\eta(t) := \sup\left\{n \geq 0 : \sum_{k=1}^n \sigma_n / r(Z_n^y) \leq t\right\} \qquad (t \geq 0). \tag{3.71}$$

Then $\eta(t) \leq n_\infty$ for all $t \geq 0$ and

$$Y_t^y = Z_{\eta(t)}^y \qquad (t \geq 0). \tag{3.72}$$

Conditioning on the embedded chain $Z^y = (Z_n^y)_{0 \leq n \leq n_\infty}$ yields

$$\begin{aligned}
\mathbb{P}[0 < |Y_t^y| < N] &= \sum_{n=0}^\infty \mathbb{P}[0 < |Z_n^y| < N, \ \eta(t) = n] \\
&= \int \mathbb{P}[Z^y \in dz] \sum_{n : 0 < |z_n| < N} \mathbb{P}[\eta(t) = n \mid Z^y = z].
\end{aligned} \tag{3.73}$$

Hence, in order to prove (3.69), it suffices to show for each $z = (z_n)_{0 \leq n \leq n_\infty}$,

$$\begin{aligned}
\lim_{t \to \infty} \sup_{s \in [0,S]} \sum_{n : 0 < |z_n| < N} &|\mathbb{P}[\eta(t) = n \mid Z^y = z] \\
&- \mathbb{P}[\eta(t+s) = n \mid Z^y = z]| = 0.
\end{aligned} \tag{3.74}$$

If $0 < |z_n| < N$ for finitely many $n$, then

$$
\sum_{n:\, 0 < |z_n| < N} |\mathbb{P}[\eta(t) = n \mid Z^y = z] - \mathbb{P}[\eta(t+s) = n \mid Z^y = z]|
$$
(3.75)
$$
\leq \mathbb{P}[0 < z_{\eta(t)} < N \mid Z^y = z] + \mathbb{P}[0 < z_{\eta(t+s)} < N \mid Z^y = z],
$$

which tends to zero as $t \to \infty$, uniformly for all $s \in [0, S]$. If $0 < |z_n| < N$ for infinitely many $n$, then we estimate

$$
\sum_{n:\, 0 < |z_n| < N} |\mathbb{P}[\eta(t) = n \mid Z^y = z] - \mathbb{P}[\eta(t+s) = n \mid Z^y = z]|
$$
(3.76)
$$
\leq \sum_{n=0}^{\infty} |\mathbb{P}[\eta(t) = n \mid Z^y = z] - \mathbb{P}[\eta(t+s) = n \mid Z^y = z]|
$$
$$
= \|\mathbb{P}[\eta(t) \in \cdot \mid Z^y = z] - \mathbb{P}[\eta(t+s) \in \cdot \mid Z^y = z]\|,
$$

where $\|\cdot\|$ denotes the total variation norm. Since $0 < |z_n| < N$ for infinitely many $n$, it is easy to see that $\liminf_{n \geq 0} r(z_n) < \infty$. Therefore, by Lemma 15 below,

$$
\text{(3.77)} \quad \lim_{t \to \infty} \sup_{s \in [0, S]} \|\mathbb{P}[\eta(t) \in \cdot \mid Z^y = z] - \mathbb{P}[\eta(t+s) \in \cdot \mid Z^y = z]\| = 0,
$$

as required. $\square$

LEMMA 15 (Coupling of exponential variables). *Let $(\sigma_n)_{n \geq 1}$ be independent, exponentially distributed random variables with mean one and let $(\lambda_n)_{n \geq 1}$ be nonnegative constants. Set*

$$
\text{(3.78)} \qquad \eta(t) := \sup\left\{ n \geq 0 : \sum_{k=1}^{n} \lambda_k \sigma_k \leq t \right\}.
$$

*If $\sum_{k=1}^{\infty} \lambda_k^2 = \infty$, then*

$$
\text{(3.79)} \qquad \lim_{t \to \infty} \sup_{s \in [0, S]} \|\mathbb{P}[\eta(t) \in \cdot] - \mathbb{P}[\eta(t+s) \in \cdot]\| = 0,
$$

*where $\|\cdot\|$ denotes the total variation norm.*

PROOF. Our lemma is very similar to [20], Lemma 2.2, except that the latter uses the additional technical assumption that $\sup_{k \geq 1} \lambda_k/\lambda_{k+1} < \infty$. Since we do not want to assume this, our proof will be a bit different from the proof there.

Let $(\eta(t))_{t \geq 0}$ and $(\eta'(t))_{t \geq s}$ be continuous-time Markov processes on $\mathbb{N}$ that jump from $k-1$ to $k$ with rate $1/\lambda_k$. We start these processes at time $0$ and time $-s$ in $\eta(0) = 0$ and $\eta'(-s) = 0$, respectively. From time $0$ onward, we let them run independently until the first time they meet, after which they are equal. We claim that this coupling is succesful. To see this, set

$$\Delta_t := \sum_{k=\eta(t)+1}^{\eta'(t)} \lambda_k \qquad (t \geq 0),$$

(3.80)
$$\tau_0 := \inf\{t \geq 0 : \Delta_t = 0\},$$

$$\tau_R := \inf\{t \geq 0 : \Delta_t \geq R\} \qquad (R > 0).$$

Then $(\Delta_t)_{t \geq 0}$ is a nonnegative, square integrable martingale. Moreover, the process $M_t := \Delta_t^2 - \langle \Delta \rangle_t$ is a martingale, where

(3.81)
$$\langle \Delta \rangle_t := \int_0^{t \wedge \tau_0} (\lambda_{\eta(u)+1} + \lambda_{\eta'(u)+1}) \, du.$$

By optional stopping, it follows that

(3.82) $\qquad \mathbb{E}[\langle \Delta \rangle_{t \wedge \tau_R}] = \mathbb{E}[\Delta_{t \wedge \tau_R}^2] \leq \mathbb{E}[(R \vee \eta'(0))^2] \qquad (t \geq 0),$

so letting $t \uparrow \infty$, we see that $\langle \Delta \rangle_{\tau_R} < \infty$ a.s. for each $R > 0$. Since $(\Delta_t)_{t \geq 0}$ is a nonnegative martingale, it has an a.s. limit as $t \to \infty$, hence, a.s. $\tau_R = \infty$ for some finite random $R > 0$. By (3.81), it follows that

(3.83) $\qquad \int_0^{\tau_0} (\lambda_{\eta(u)+1} + \lambda_{\eta'(u)+1}) \, du < \infty \qquad \text{a.s.}$

Since $\sum_{k=1}^\infty \lambda_k^2 = \infty$, it is not hard to see that $\int_0^\infty \lambda_{\eta(u)+1} \, du = \infty$ a.s. Therefore, (3.83) implies that $\tau_0 < \infty$ a.s.

This shows that our coupling is successful, hence, the total variation distance between $\mathbb{P}[\eta(t) \in \cdot]$ and $\mathbb{P}[\eta'(t) \in \cdot]$ tends o zero as $t \to \infty$. It is easy to see that the stopping time $\tau_0$ is stochastically increasing in $s$, hence, our estimates are uniform for all $s \in [0, S]$, for each $S > 0$. Since $\eta(t)$ is distributed as the process in (3.78) and $\eta'(t)$ is distributed as $\eta(t+s)$, our claim follows. □

PROOF OF THEOREM 12. It suffices to show that

(3.84) $$\lim_{n \to \infty} \mathbb{P}[0 < |Y_{t_n}^y| < N] = 0$$

for any sequence of times $t_n \to \infty$. By Lemma 14 and a diagonal argument, we can choose $T_n \to \infty$ such that

(3.85) $\qquad \lim_{n \to \infty} \left| \mathbb{P}[0 < |Y_{t_n}^y| < N] - \frac{1}{T_n} \int_0^{T_n} ds \, \mathbb{P}[0 < |Y_{t_n+s}^y| < N] \right| = 0.$

The proof of Proposition 13 actually works more generally than for Césaro times only. Let $\sigma_{T_n}$ and $\sigma_S$ be uniformly distributed on $[0, T_n]$ and $[0, S]$, respectively, independent of each other and of $Y^y$. Since for each fixed $S > 0$, the total variation distance between $\mathbb{P}[t_n + \sigma_{T_n} \in \cdot]$ and $\mathbb{P}[t_n + \sigma_{T_n} + \sigma_S \in \cdot]$ tends to zero as $n \to \infty$, the proof of Proposition 13 tells us that

$$\lim_{n \to \infty} \mathbb{P}[0 < |Y^y_{t_n + \sigma_{T_n}}| < N] = 0. \tag{3.86}$$

Combining this with (3.85), we arrive at (3.55). □

3.5. *Homogeneous invariant laws.* In this section we prove Theorem 3, using Theorem 6, Corollary 11 and Theorem 12. Since Theorem 6 speaks about nontrivial laws, while we are interested in coexisting and nonzero laws, we need the following lemma. Recall the definitions of local nonsingularity and $X$-nontriviality above Theorem 6.

LEMMA 16 (Local nonsingularity). *Let $X$ be either the neutral Neuhauser–Pacala model, affine voter model or rebellious voter model, and let $Y$ be its dual parity preserving branching process. Then, for any $0 \leq \alpha \leq 1$:*

  (a) *Each homogeneous coexisting law on $\{0,1\}^{\mathbb{Z}^d}$ is $X$-nontrivial.*
  (b) *Each homogeneous nonzero law on $\{0,1\}^{\mathbb{Z}^d}$ is $Y$-nontrivial.*

Lemma 16 can be verified by simple, but lengthy considerations, which we leave to the reader. To show that this point requires some care, we warn the reader that the claim in (a) does not hold for the disagreement voter model with $\alpha = 0$, since for this model, the alternating configurations $\ldots 0101010 \ldots$ are traps.

PROOF OF THEOREM 3. We start with part (a). Since $Y$ survives, by Lemma 1, $\nu_X^{1/2}$ is not concentrated on $\{\underline{0}, \underline{1}\}$. Conditioning on the complement of the set $\{\underline{0}, \underline{1}\}$, we obtain a homogeneous coexisting invariant law of $X$. By Lemma 16, each homogeneous coexisting invariant law of $X$ is $X$-nontrivial in the sense of Theorem 6. Therefore, since $\alpha < 1$, it follows from Theorem 6 and Corollary 11(a) that $\nu_X^{1/2}$ is the only homogeneous coexisting invariant law. If, moreover $\alpha > 0$ and $Y$ is not stable, then by Theorem 6 and Theorem 12, the convergence in (1.13) holds. Here, condition (3.54) of Theorem 12 is easily seen to follow from our assumption that $\alpha > 0$.

To prove also part (b), we observe that since $X$ survives, by Lemma 1, $\nu_Y^{1/2}$ is not concentrated on $\underline{0}$, hence, $\nu_Y^{1/2}$ conditioned on being nonzero is a homogeneous nonzero invariant law of $Y$. By Lemma 16, each homogeneous nonzero invariant law of $Y$ is $Y$-nontrivial. Therefore, since $\alpha < 1$, it follows from Theorem 6 and Corollary 11(a) that $\nu_Y^{1/2}$ is the only homogeneous nonzero invariant law. If, moreover, $\alpha > 0$ and $d \geq 2$, then by Theorem 6 and Corollary 11(b), the convergence in (1.14) holds. □

3.6. *A counterexample to a result by Simonelli.* Our Theorem 6 is similar to Theorem 1 in [22]. The proof is also similar, with his Theorem 2 playing the same role as our Corollary 9. An important difference is that while we use the "norm" $\|\cdot\|_{\mathcal{B}}$ defined in (3.11), Simonelli works with the usual $\ell^1$-norm $|\cdot|$. As we will see in a moment, the result of this is that his Theorem 2 is false. We do not know if his main result Theorem 1 is correct or not; as it is, his proof depends on his Theorem 2, and is therefore not correct.

More precisely, in [22] it is assumed that $X$ and $Y$ are interacting particle systems on $\mathbb{Z}^d$ satisfying a duality relation of the form (1.8). It is assumed that

$$(3.87) \qquad 0 < \mathbb{P}[Y_t^y(i) = 1] < 1 \qquad (t > 0, i \in \mathbb{Z}^d, y \neq \underline{0}),$$

where $Y^y$ denotes the process $y$ started in $Y_0^y = y$. Then, translated into our terminology, Simonelli ([22], Theorem 2) states that:

> If $\mathbb{P}[Y_0 \in \cdot]$ is homogeneous and nonzero, and $x_n \in \{0,1\}^{\mathbb{Z}^d}$ satisfy $\lim_{n\to\infty} |x_n| = \infty$, then $\lim_{n\to\infty} \mathbb{P}[|Y_t x_n| \text{ is odd}] = \frac{1}{2}$ for all $t > 0$.

This claim is false, as is shown by the following:

COUNTEREXAMPLE. We pick some $\alpha \in (0,1)$ and we take for $Y$ the ADBARW, which is the dual of the rebellious voter model. It is easy to see that $Y$ satisfies (3.87). Set $\mathbb{P}[Y_0 = \underline{1}] = 1$, which is homogeneous and nonzero, and $x_n(i) := 1$ for $i = 1, \ldots, n$ and $x_n(i) := 0$ otherwise. Using duality, we see that $\mathbb{P}[|Y_t x_n| \text{ is odd}] = \mathbb{P}[|Y_0 X_t^{x_n}| \text{ is odd}]$, where $X^{x_n}$ denotes the rebellious voter model started in $X_0^{x_n} = x_n$. The dynamics of the rebellious voter model are such that at $t = 0$, the only places where sites can flip are at the endpoints of the interval $\{1, \ldots, n\}$. Therefore, it is easy to see that $\mathbb{P}[X_t^{x_n} = x_n] \geq 1 - e^{-4t}$, and as a result,

$$(3.88) \qquad \lim_{t\to 0} \mathbb{P}[|Y_0 X_t^{x_n}| \text{ is odd}] = \begin{cases} 1, & \text{if } n \text{ is odd}, \\ 0, & \text{if } n \text{ is even}, \end{cases}$$

where the convergence is *uniform in $n$*. It follows that, for $t$ sufficiently small, the limit $\lim_{n\to\infty} \mathbb{P}[|Y_t x_n| \text{ is odd}]$ does not exist.

There seems to be no easy way to repair Simonelli's Theorem 2. The essential observation behind our Lemma 7 is that, instead of $x_n$ being large, one needs that $x_n$ contains many places where parity can change. For the rebellious voter model, this means that $x_n$ must contain many places where the two types meet. Our proof of Lemma 7 differs substantially from the methods used in [22].

The place where Simonelli's proof goes wrong is his formula (5), where it is claimed that if $\mathbb{P}[Y_0 \in \cdot]$ is homogeneous and nonzero, then for each $\varepsilon > 0$ and $t > 0$ there exists a $\delta > 0$ such that, for all $s \in [0, \varepsilon]$,

$$\mathbb{P}[Y_t \in \{y : \delta < \mathbb{P}[Y_s^y(0) = 1] < 1 - \delta\}] > 1 - \varepsilon \qquad ([22], (5)).$$

While this inequality is true for any $s > 0$ fixed with $\delta = \delta(t, s, \varepsilon)$ depending on $s, t$, and $\varepsilon$ due to his assumption (3.87) and continuity of probability measures, $\delta$ cannot be chosen uniformly in $s \in [0, \varepsilon]$. In fact, for $s = 0$, the set $\{y : \delta < \mathbb{P}[Y_s^y(0) = 1] < 1 - \delta\}$ is empty for any $\delta > 0$.

**4. Complete convergence of the rebellious voter model.** In this section we prove coexistence and complete convergence for the rebellious voter model, as stated in Theorem 4. The main tool for this will be Theorem 5, which states that the dual ADBARW dominates oriented percolation. This theorem is formulated and proven in the following Section 4.1.

4.1. *Comparison with oriented percolation.* We let $Y$ be an ADBARW started in an arbitrary deterministic initial state $Y_0 = y \in \{0, 1\}^{\mathbb{Z}}$, and we define sets of "good" points $(\chi_n)_{n \geq 0}$ as in (2.13). We start by considering a single time step in the case that $\alpha = 0$. Our first result says that "good" events have a large probability.

PROPOSITION 17 (Good events are probable). *Assume that $\alpha = 0$. Then, for each $p < 1$, there exist $L \geq 1$ and $T > 0$, not depending on the initial state $y$, such that $\mathbb{P}[x \in \chi_1] \geq p$ for all $x \in \mathbb{Z}_{\text{odd}}$ such that $\chi_0 \cap \{x - 1, x + 1\} \neq \varnothing$.*

PROOF. Our basic observation is that, in case $\alpha = 0$, the right-most particle of an ADBARW started in a finite initial state has a drift to the right. Indeed, if $r_t := \max\{i \in \mathbb{Z} : Y_t(i) = 1\}$, then depending on the configuration near the right-most particle, the changes in $r_t$ due to the various possible jumps and the resulting drift are as follows:

(4.1)

| configuration | change in $r_t$ | rate | drift |
|---|---|---|---|
| ...00100... | $+2$ | $1$ | $+2$ |
| ...01100... | $\begin{cases} +2 \\ +1 \end{cases}$ | $\begin{cases} 1 \\ 1 \end{cases}$ | $+3$ |
| ...10100... | $\begin{cases} +2 \\ -1 \end{cases}$ | $\begin{cases} 1 \\ 1 \end{cases}$ | $+1$ |
| ...11100... | $\begin{cases} +2 \\ +1 \\ -2 \end{cases}$ | $\begin{cases} 1 \\ 1 \\ 1 \end{cases}$ | $+1.$ |

This shows that in each configuration the drift is at least one. We note, however, that it is not possible to stochastically bound $r_t$ from the left by a random walk with positive drift (independent of anything else). This will cause a slight complication in what follows; in fact, we will use two random walks that become active when $Y(r_t - 1) = 0$ or 1, respectively [see formula (4.11) below].

We need to prove that if $\chi_0 \cap \{x - 1, x + 1\} \neq \varnothing$ for some $x \in \mathbb{Z}_{\text{odd}}$, then $\mathbb{P}[x \in \chi_1] \geq p$. By symmetry, we may without loss of generality assume that $x = 1$ and $0 \in \chi_0$. To simplify notation, let us identify subsets of $\mathbb{Z}$ with their indicator functions. Then, assuming that $y \cap \{-L, \ldots, L\} \neq \varnothing$, which is equivalent to $0 \in \chi_0$, we need to show that the probability

(4.2) $\quad \mathbb{P}[Y_T \cap \{L, \ldots, 3L\} \neq \varnothing$ and $Y_t \cap \{-2L, \ldots, 6L\} \neq \varnothing \ \forall 0 < t < T]$

can be made arbitrarily large by choosing $L$ and $T$ appropriately. In view of this, we are actually not interested in the right-most particle of our AD-BARW, but in the particle that is closest to our target $2L$. Thus, we put

(4.3) $\qquad s_t := \inf\{i \geq 0 : Y_t(2L - i) \vee Y_t(2L + i) = 1\}.$

Assuming that $s_0 \leq 3L$ which follows from $y \cap \{-L, \ldots, L\} \neq \varnothing$, we need to show that the probability

(4.4) $\qquad \mathbb{P}\left[s_T \leq L, \sup_{0 \leq t \leq T} s_t \leq 4L\right]$

can be made arbitrarily large. For any $n \geq 0$, we set

(4.5) $\qquad \begin{aligned} \tau_{\leq n} &:= \inf\{t \geq 0 : s_t \leq n\}, \\ \tau_{\geq n} &:= \inf\{t \geq 0 : s_t \geq n\}. \end{aligned}$

We observe from (4.1) that whenever an ADBARW borders at least two empty sites, it tends to invade these with a drift of at least one. In view of this, we choose $T = 2L$. By Lemmas 18 and 19 below, there exist constants $C, \lambda > 0$ such that

(4.6) $\qquad \mathbb{P}[\tau_{\leq 2} \leq T \text{ and } \tau_{\geq 4L} \geq T] \geq 1 - Ce^{-\lambda L}.$

Using Lemma 19 once more, we see that, moreover, for some $C', \lambda' > 0$,

(4.7) $\qquad \mathbb{P}[s_t \leq L \ \forall \tau_{\leq 2} \leq t \leq T] \geq 1 - C'e^{-\lambda' L}.$

Combining these two estimates, we see that the probability in (4.4) can be made as close to one as one wishes by choosing $L$ large enough. $\square$

We still need to prove two lemmas.

LEMMA 18 (Hitting the target). *For each $\delta > 0$, there exist constants $C, \lambda > 0$ such that if $s_0 \leq K$, then*

(4.8) $\qquad \mathbb{P}[\tau_{\leq 2} \geq (1 + \delta)K] \leq Ce^{-\lambda K} \qquad (K \geq 1).$

PROOF. Let

(4.9) $\qquad s_t^- := \inf\{i \geq 0 : Y_t(2L - i) = 1\}.$

Since the evolution rules of the process are symmetric, we may assume without loss of generality that $s_0^- \leq K$. We set

$$(4.10) \qquad \phi_q(t) := \int_0^t 1_{\{Y(s_t^- - 1) = q\}} \, du \qquad (t \geq 0, q = 0, 1).$$

Let $R^0, R^1$ be continuous-time random walks on $\mathbb{Z}$, starting in zero, with the following jump rates:

$$(4.11) \qquad \begin{array}{ccc} \text{random walk} & \text{jump size} & \text{rate} \\ R_t^0 & \begin{cases} +2 \\ -1 \end{cases} & \begin{array}{c} 1 \\ 1 \end{array} \\ R_t^1 & \begin{cases} +2 \\ +1 \\ -2 \end{cases} & \begin{array}{c} 1 \\ 1 \\ 1. \end{array} \end{array}$$

In view of (4.1), we can couple $s_t^-$ to $R^0, R^1$ in such a way that

$$(4.12) \qquad s_t^- \leq s_0^- - R_{\phi_0(t)}^0 - R_{\phi_1(t)}^1 \qquad (0 \leq t \leq \tau_{\leq 2}),$$

where $R^0$ and $R^1$ are independent of each other and of $\phi_0, \phi_1$. It follows from large deviation theory, more precisely, from Cramér's theorem (Theorem 27.3 in [16]) and a little calculation, that for each $\varepsilon > 0$ there exist constants $C_\varepsilon$ and $\lambda_\varepsilon > 0$ such that

$$(4.13) \qquad \mathbb{P}[|R_t^q - t| \geq \varepsilon t] \leq C_\varepsilon e^{-\lambda_\varepsilon t} \qquad (t \geq 0, q = 0, 1).$$

We claim that there exist $C'_\varepsilon$ and $\lambda'_\varepsilon$ such that

$$(4.14) \qquad \mathbb{P}[|R_s^q - s| \geq \varepsilon t] \leq C'_\varepsilon e^{-\lambda'_\varepsilon t} \qquad (0 \leq s \leq t, q = 0, 1).$$

Since it is not obvious that $\mathbb{P}[|R_s^q - s| \geq \varepsilon t] \leq \mathbb{P}[|R_t^q - t| \geq \varepsilon t]$, it is not entirely trivial to get from (4.13) to (4.14). Here is a clumsy argument: If $\frac{1}{2}t \leq s \leq t$, then

$$(4.15) \quad \mathbb{P}[|R_s^q - s| \geq \varepsilon t] \leq \mathbb{P}[|R_s^q - s| \geq \varepsilon s] \leq C_\varepsilon e^{-\lambda_\varepsilon s} \leq C_\varepsilon e^{-(1/2)\lambda_\varepsilon t}.$$

If $0 \leq s \leq \frac{1}{2}t$, then, by the independence of random walk increments and what we have just proved,

$$(4.16) \qquad \begin{aligned} \mathbb{P}[|R_t^q - t| \geq \tfrac{1}{2}\varepsilon t] &\geq \mathbb{P}[|R_s^q - s| \geq \varepsilon t] \mathbb{P}[|R_{t-s}^q - (t-s)| \leq \tfrac{1}{2}\varepsilon t] \\ &\geq \mathbb{P}[|R_s^q - s| \geq \varepsilon t](1 - C_{\varepsilon/2} e^{-(1/2)\lambda_{\varepsilon/2} t}), \end{aligned}$$

hence, by (4.13),

$$(4.17) \qquad \mathbb{P}[|R_s^q - s| \geq \varepsilon t] \leq (1 - C_{\varepsilon/2} e^{-(1/2)\lambda_{\varepsilon/2} t})^{-1} C_\varepsilon e^{-\lambda_{\varepsilon/2} t}.$$

Combining (4.15) and (4.17), we obtain (4.14).

To prove (4.8), we set $M := (1+\delta)K$, we choose $\varepsilon$ such that $(1-\varepsilon)M = K$, and observe that by (4.12), the fact that $\phi_0(M) + \phi_1(M) = M$, and (4.14),

$$
\begin{aligned}
\mathbb{P}[\tau_{\leq 2} \geq M] &\leq \mathbb{P}[R^0_{\phi_0(M)} + R^1_{\phi_1(M)} \leq (1-\varepsilon)M] \\
&\leq \mathbb{P}[|R^0_{\phi_0(M)} + R^1_{\phi_1(M)} - \phi_0(M) - \phi_1(M)| \geq \varepsilon M] \\
&\leq \sum_{q=0}^{1} \mathbb{P}[|R^q_{\phi_q(M)} - \phi_q(M)| \geq \tfrac{1}{2}\varepsilon M] \leq 2C'_{\varepsilon/2} e^{-\lambda'_{\varepsilon/2} M}.
\end{aligned}
\tag{4.18}
$$

Setting $C := 2C'_{\varepsilon/2}$ and $\lambda := \lambda'_{\varepsilon/2}(1+\delta)$, we arrive at (4.8). $\square$

LEMMA 19 (Escaping the target). *For each $\delta > 0$, there exist constants $C, \lambda > 0$ such that if $s_0 \leq K$, then*

$$\mathbb{P}[\tau_{\geq(1+\delta)K} \leq T] \leq C(T+1)e^{-\lambda K} \qquad (T > 0,\ K \geq 1). \tag{4.19}$$

PROOF. We start by showing that we can choose $\lambda > 0$ such that if $s_0 \leq L' \leq L$, then

$$\mathbb{P}[\tau_{\leq 2} > \tau_{\geq L}] \leq e^{-\lambda(L-L')}. \tag{4.20}$$

Indeed, it is not hard to see from (4.1) that, for $\lambda$ sufficiently small, the process

$$M_t := e^{\lambda s_{\inf\{t, \tau_{\leq 2}\}}} \qquad (t \geq 0) \tag{4.21}$$

is a supermartingale. Setting $\tau := \inf\{\tau_{\leq 2}, \tau_{\geq L}\}$, by optional stopping, it follows that

$$e^{\lambda L} \mathbb{P}[s_\tau \geq L] \leq \mathbb{E}[M_\tau] \leq e^{\lambda L'}, \tag{4.22}$$

which proves (4.20). To prove (4.19), we note that each time the process $s_t$ enters $\{0, 1, 2\}$, it must stay there at least an exponentially distributed time with mean $\tfrac{1}{2}$. Therefore, the number of excursions from $\{0, 1, 2\}$ during a time interval of length $T$ is bounded by a Poisson random variable with mean $2T$, and by (4.20), the number of excursions from $\{0, 1, 2\}$ that go beyond $(1+\delta)K$ is bounded by a Poisson random variable $W$ with

$$\mathbb{P}[W > 0] \leq \mathbb{E}[W] = 2T e^{-\lambda((1+\delta)K - 2)}. \tag{4.23}$$

Using (4.19) once more, we conclude that the probability in (4.19) is bounded by $e^{-\lambda \delta K} + 2T e^{-\lambda((1+\delta)K - 2)}$. $\square$

We now turn to the proof of Theorem 5. It suffices to prove the statement for a single time step; the general statement then follows by induction. Thus, we need to show that for each $p < 1$ we can find $L, T$ so that we can define i.i.d. Bernoulli random variables $\{\omega_{(x,1)} : x \in \mathbb{Z}_{\mathrm{odd}}\}$ with $\mathbb{P}[\omega_{(x,1)} = 1] = p$ such that $x \in \chi_1$ whenever $\omega_{x,1} = 1$ and $\chi_0 \cap \{x-1, x+1\} \neq \varnothing$. Since we can thin the set of points where $\omega_{(x,1)} = 1$ if necessary, it suffices if the $\{\omega_{(x,1)} : x \in \mathbb{Z}_{\text{odd}}\}$ are independent and $\mathbb{P}[\omega_{(x,1)} = 1] \geq p$ for each $x$. It actually suffices if the $\{\omega_{(x,1)} : x \in \mathbb{Z}_{\text{odd}}\}$ are $m$-dependent for some fixed $m \geq 1$, since a well-known result (see [19], Theorem B26) tells us that $m$-dependent random variables with intensity $p$ can be estimated from below by independent random variables with intensity $p' = p'(p, m)$ depending on $p$ and $m$ in such a way that $\lim_{p \to 1} p'(p, m) = 1$. See also [12] who considers a more general form of $m$-dependence for oriented percolation and the comparison argument.

PROOF OF THEOREM 5. As mentioned before, it suffices to prove the statement for a single time step, and for the process started in any deterministic initial state $Y_0 = y$. A naive approach is to put $\omega_{(x,1)} := 1_{\{x \in \chi_1\}}$ for $x \in J := \{x \in \mathbb{Z} : \chi_0 \cap \{x-1, x+1\} \neq \varnothing\}$ and $\omega_{(x,1)} := 1$ otherwise. However, since these events are not $m$-dependent for any fixed $m$, we are going to extend our definition of a "good" event, in such a way that the new events still have a high probability, and are $m$-dependent. To that aim, we put

$$\chi'_1 := \{x \in \chi_1 : \text{there is no path in the graphical representation ending}$$
(4.24) $$\text{at time } T \text{ in } \{2Lx - 4L, \ldots, 2Lx + 4L\} \text{ and starting}$$
$$\text{at time } 0 \text{ outside } \{2Lx - 11L, \ldots, 2Lx + 11L\}\},$$

where, as before, we choose $T = 2L$. The motivation for this is as follows. Let

(4.25) $$r_t := \sup\{i \in \mathbb{Z} : \text{there is a path starting at } (i, T-t)$$
$$\text{and ending in } \{T\} \times \{2Lx - 4L, \ldots, 2Lx + 4L\}\}.$$

It is not hard to see that $r_t$ can be bounded from above by a random walk that makes jumps of size $+1$ and $+2$, both with rate 1. Therefore, the expected distance covered by such a random walk is $3T = 6L$. A large deviation estimate of the same sort as used in the proof of Lemma 18 now tells us that for each $x \in \mathbb{Z}_{\text{odd}}$, the probability of there being a path ending at time $T$ in $\{2Lx - 4L, \ldots, 2Lx + 4L\}$ and starting at time 0 outside $\{2Lx - 11L, \ldots, 2Lx + 11L\}$ tends to zero exponentially fast as $L \to \infty$. Combining this with Proposition 17, we see that for each $p$ we can choose $L, T$ in such a way that $\mathbb{P}[x \in \chi'_1] \geq p$ for all $x \in J$. Choosing $\alpha'$ close enough to zero so that the probability of any event with rate $\alpha'$ happening in the graphical representation in the block $[0, T] \times \{2Lx - 11L, \ldots, 2Lx + 11L\}$ is small, we conclude that

(4.26) For each $p < 1$ there exist an $\alpha' > 0$ such that for all $\alpha \in [0, \alpha')$ there exist $L \geq 1$ and $T > 0$ such that $\mathbb{P}[x \in \chi'_1] \geq p$ for all $x \in J$.

We now put $\omega_{(x,1)} := 1_{\{x \in \chi'_1\}}$ for $x \in J$ and $\omega_{(x,1)} := 1$ otherwise. The $\{\omega_{(x,1)} : x \in \mathbb{Z}_{\text{odd}}\}$ constructed in this way are $m$-dependent for a suitable $m$ [in fact, $\omega_{(x,1)}$ and $\omega_{(x',1)}$ are independent if $|x - x'| \geq 12$], so by the arguments preceding this proof, they can be bounded from below by i.i.d. Bernoulli random variables with an intensity that can be made arbitrarily high. $\square$

4.2. *Complete convergence.* In this section we prove Theorem 4 about coexistence and complete convergence for the rebellious voter model $X$ for small enough $\alpha$. Throughout this section we choose some $p \in (p_c, 1)$, where $p_c$ is the critical value for survival of oriented percolation, and we fix $\alpha' > 0$, $L \geq 1$, and $T > 0$ such that the process $(\chi_n)_{n \geq 0}$ defined in (2.13) can be coupled to an oriented percolation process $(W_n)_{n \geq 0}$ with parameter $p$ as in Theorem 5. The proof of Theorem 4 is based on the following two lemmas. Recall from Section 2.1 that the ADBARW is both the dual and the interface model of the rebellious voter model.

LEMMA 20 (Almost sure extinction versus unbounded growth). *For all $\alpha \in [0, \alpha')$, the ADBARW $Y^y$ started in any finite initial state $Y_0^y = y$ satisfies*

$$(4.27) \quad \mathbb{P}\left[\lim_{t \to \infty} |Y_t^y| = \infty \text{ or } \exists t \geq 0 \text{ s.t. } Y_t^y = \underline{0}\right] = 1 \quad (|y| < \infty).$$

*Moreover,* $\mathbb{P}[\lim_{t \to \infty} |Y_t^y| = \infty] > 0$ *for all* $y \neq \underline{0}$.

LEMMA 21 (Intersection of independent processes). *For all $\alpha \in [0, \alpha')$, two independent ADBARW's $Y$ and $Y'$ satisfy for $t \geq 0, N \geq 1$*

$$(4.28) \quad \lim_{t \to \infty} \mathbb{P}[|Y_t Y_t'| \geq N] = \mathbb{P}[Y_t \neq \underline{0} \ \forall t \geq 0] \cdot \mathbb{P}[Y_t' \neq \underline{0} \ \forall t \geq 0].$$

PROOF OF THEOREM 4. If $\alpha$ is as in Lemma 20, then the ADBARW is unstable and survives. Therefore, by Lemmas 1 and 2, the rebellious voter model exhibits coexistence, survival, and its dual ADBARW is not stable. This proves part (a).

To prove part (b), it suffices to consider deterministic initial states. Let $X^x$ denote the rebellious voter model started in $X_0^x = x$ and set $\rho_q(x) := \mathbb{P}[X_t^x = \underline{q} \text{ for some } t \geq 0] \ (q = 0, 1)$. By (2.4), it suffices to show that, for each $|y| < \infty$,

$$(4.29) \quad \lim_{t \to \infty} \mathbb{P}[|X_t^x y| \text{ is odd }] \\ = (1 - \rho_0(x) - \rho_1(x))\tfrac{1}{2}\mathbb{P}[Y_t^y \neq \underline{0} \ \forall t \geq 0] + \rho_1(x)1_{\{|y| \text{ is odd}\}}.$$

Since the dynamics of $X$ are symmetric for the 0's and 1's, we may without loss of generality assume that the starting configuration $x$ has infinitely many 0's. In this case $\rho_1(x) = 0$ and (4.29) reduces to

$$(4.30) \quad \lim_{t \to \infty} \mathbb{P}[|X_t^x y| \text{ is odd }] = \tfrac{1}{2}(1 - \rho_0(x))\mathbb{P}[Y_t^y \neq \underline{0} \ \forall t \geq 0].$$

Since, by duality,

$$(4.31) \quad \begin{aligned} &\mathbb{P}[|X_t^x y| \text{ is odd }] \\ &= \int \mathbb{P}[X_{(t-1)/2}^x \in d\tilde{x}] \int \mathbb{P}[Y_{(t-1)/2}^y \in d\tilde{y}] \mathbb{P}[|X_1^{\tilde{x}} Y_0^{\tilde{y}}| \text{ is odd}], \end{aligned}$$

by Lemma 7, in order to prove (4.30), it suffices to show that

$$(4.32) \quad \begin{aligned} \lim_{t \to \infty} \mathbb{P}[|\{i \in \mathbb{Z} : Y_{(t-1)/2}^y(i) = 1, \\ X_{(t-1)/2}^x(i+1) \neq X_{(t-1)/2}^x(i+2)\}| \geq N] \\ = \mathbb{P}[X_t^x \neq \underline{0} \ \forall t \geq 0] \cdot \mathbb{P}[Y_t^y \neq \underline{0} \ \forall t \geq 0] \quad (N \geq 1). \end{aligned}$$

Let $Y_t'(i) := 1_{\{X_t^x(i+1) \neq X_t^x(i+2)\}}$, that is, $Y'$ is the interface model of $X^x$, translated over a distance one. Then (4.32) simplifies for $N \geq 1$ to

$$(4.33) \quad \lim_{t \to \infty} \mathbb{P}[|Y_{(t-1)/2}^y Y_{(t-1)/2}'| \geq N] = \mathbb{P}[Y_t' \neq \underline{0} \ \forall t \geq 0] \cdot \mathbb{P}[Y_t^y \neq \underline{0} \ \forall t \geq 0],$$

which is true by Lemma 21. □

The rest of this section is occupied by the proofs of Lemmas 20 and 21.

PROOF OF LEMMA 20. By Lemma 10(b), the claims will follow provided that we show that

$$(4.34) \quad \inf_{|y|>0} \mathbb{P}\left[\lim_{t \to \infty} |Y_t^y| = \infty\right] > 0.$$

By translation invariance, it suffices to consider the infimum over all $y$ such that $y(0) = 1$. Choose $p > p_c$, the critical value for survival of oriented percolation. By Theorem 5, the process $(\chi_n)_{n \geq 0}$ defined in (2.13) can be coupled to an oriented percolation process $(W_n)_{n \geq 0}$ such that $W_0 = \chi_0$ and $W_n \subset \chi_n$ for all $n \geq 1$. Since $0 \in \chi_0$ due to the fact that $y(0) = 1$, and since $p > p_c$, there is a positive probability that $\lim_{n \to \infty} |W_n| = \infty$ and, hence, $\lim_{t \to \infty} |Y_t^y| = \infty$. □

The proof of Lemma 21 is somewhat more involved. We start with some preparatory lemmas. Our first lemma says that if $y$ and $y'$ are close in many places, then $|Y_t^y y'|$ is large with probability close to one.

LEMMA 22 (Charging target sets). *Assume that $\alpha < 1$, let $Y^y$ be an AD-BARW started in $y$, let $y' \in \{0,1\}^{\mathbb{Z}}$, and let $K \geq 1$, $t > 0$. Set $D_K(y,y') := \{(i,j) \in \mathbb{Z}^2 : y(i) = 1 = y'(j), |i-j| \leq K\}$. Then*

$$(4.35) \qquad \lim_{M \to \infty} \inf_{|D_K(y,y')| \geq M} \mathbb{P}[|Y_t^y y'| \geq N] = 1 \qquad (N \geq 1).$$

PROOF. Set $C := \{j : \exists i \text{ s.t. } (i,j) \in D_K(y,y')\}$. For each $j \in C$, choose in some unique way a site $i$ with $y(i) = 1$ for which $|i - j|$ is minimal. Let $I := \{i, \ldots, j\}$ if $i \leq j$ and $I := \{j, \ldots, i\}$ if $j \leq i$. Let $G_j$ denote the event that in the graphical representation for $Y^y$ (see Section 3.1), there is an odd number of paths from $(i,0)$ to $(j,t)$, while during the time interval $[0,t]$, there are no arrows starting outside $I$ and ending in $I$. Then, for given $K$ and $t$, the probability of $G_j$ is uniformly bounded from below. To see this, by symmetry, we may assume $i \leq j$. Then the particle at $i$ may branch to the right, producing two particles at $i+1$ and $i+2$, which can again branch to the right, creating, in a finite number of steps, a particle at $j$, while with positive probability, nothing else happens in $I$. Now, if $|D_K(y,y')| \geq M$, then we can select $C' \subset C$ such that $|j - j'| \geq 2K + 1$ for each $j, j' \in C'$ with $j \neq j'$, and $|C'| \geq M/(2K+1)^2$. Then the events $G_j$ with $j \in C'$ are independent with a probability that is uniformly bounded from below, hence, if $M$ is sufficiently large, then with large probability many of these events will occur. This proves (4.35). □

In what follows, for any $x \in \mathbb{Z}$, we define $I_x$ as in (2.12), and for $y \in \{0,1\}^{\mathbb{Z}}$, we define

$$(4.36) \qquad \eta(y) := \{x \in \mathbb{Z}_{\text{even}} : \exists i \in I_x \text{ s.t. } y(i) = 1\}.$$

LEMMA 23 (Charging target intervals). *If $Y$ is an ADBARW with parameter $\alpha \in [0, \alpha')$, started in $Y_0 = y$, then*

$$(4.37) \qquad \lim_{t \to \infty} \mathbb{P}[Y_t \neq \underline{0} \text{ and } |\eta(Y_t)| < N] = 0 \qquad (N \geq 1).$$

PROOF. By Lemma 22 applied with $K = L$ and $y' = \sum_{x \in \mathbb{Z}} \delta_{2Lx}$, for each $t > 0$,

$$(4.38) \qquad \lim_{M \to \infty} \inf_{|y| \geq M} \mathbb{P}[|\eta(Y_t^y)| \geq N] = 1 \qquad (N \geq 1).$$

Now write

$$(4.39) \quad \begin{aligned} &\mathbb{P}[Y_t \neq \underline{0} \text{ and } |\eta(Y_t)| < N] \\ &= \mathbb{P}[Y_t \neq \underline{0} \text{ and } |\eta(Y_t)| < N \mid 0 < |Y_{t-1}| < M]\mathbb{P}[0 < |Y_{t-1}| < M] \\ &\quad + \mathbb{P}[Y_t \neq \underline{0} \text{ and } |\eta(Y_t)| < N \mid |Y_{t-1}| \geq M]\mathbb{P}[|Y_{t-1}| \geq M]. \end{aligned}$$

By Lemma 20, the first term on the right-hand side of (4.39) tends to zero as $t \to \infty$, while by (4.38), the limsup as $t \to \infty$ of the second term can be made arbitrarily small by choosing $M$ large enough. □

Fix Bernoulli random variables $\{\omega_z : z \in \mathbb{Z}^2_{\text{even}}\}$ with intensity $p$ as in Section 4.1, and for each $A \subset \mathbb{Z}_{\text{even}}$, let $W^A = (W^A_n)_{n \geq 0}$ denote the oriented percolation process started in $A$ defined in (2.11). Using the same Bernoulli random variables $\omega_z$, we can define a process $\overline{W} = (\overline{W}_n)_{n \in \mathbb{Z}}$ by

$$(4.40) \qquad \overline{W}_n := \{x \in \mathbb{Z} : (x,n) \in \mathbb{Z}^2_{\text{even}}, -\infty \to (x,n)\} \qquad (n \in \mathbb{Z}),$$

where $-\infty \to (x,n)$ means that there exists an infinite open path with respect to the $\omega_z$, starting at time $-\infty$ and ending at $(x,n)$. Then $\overline{W}$ is a stationary (with respect to shifts on $\mathbb{Z}^2_{\text{even}}$) oriented percolation process. We call

$$(4.41) \qquad \overline{\nu}_{\text{W}} := \mathbb{P}[\overline{W}_{2n} \in \cdot] \qquad (n \in \mathbb{Z})$$

the *upper invariant law* of $W$. It is known that, for each $K \geq 1$ and $x \in \mathbb{Z}_{\text{even}}$,

$$(4.42) \qquad \lim_{n \to \infty} \mathbb{P}[W^{\{x\}}_n \neq \varnothing \text{ and } \overline{W}_n \cap [-K,K] \not\subset W^{\{x\}}_n] = 0.$$

This follows, for example, from [17], Theorem 2.27 (see also Theorem 2.28). The statements there are for the one-dimensional nearest neighbor voter model, but the proofs apply to our setting as well. Alternatively, one may consult [10] and [2], Theorem 5, where this (as well as the much more powerful shape theorem) is proved in a multidimensional setting, for $p$ sufficiently close to one and for $p > p_c$, respectively.

LEMMA 24 (Local comparison with upper invariant law). *Let $Y$ be an ADBARW with parameter $\alpha \in [0, \alpha')$, started in $Y_0 = y$, and let $t_k \to \infty$. Then, for each $k$, we can couple $Y_{t_k}$ to a random variable $V_k$ with law $\overline{\nu}_{\text{W}}$, in such a way that*

$$(4.43) \qquad \liminf_{k \to \infty} \mathbb{P}[V_k \cap [-K,K] \subset \eta(Y_{t_k})] \geq \mathbb{P}[Y_t \neq \underline{0} \; \forall t \geq 0] \qquad (K \geq 1).$$

PROOF. By Theorem 5, for each $s \geq 0$, we can couple $Y$ to an oriented percolation process $W^s = (W^s_n)_{n \geq 0}$ started in $W^s_0 = \eta(Y_s)$, defined by Bernoulli random variables that are independent of $(Y_u)_{u \in [0,s]}$, in such a way that

$$(4.44) \qquad W^s_n \subset \eta(Y_{s+nT}) \qquad (n \geq 0).$$

We couple $W^s$ to an "upper invariant" oriented percolation process $\overline{W}^s$ as in (4.40), and put $V_{s,n} := \overline{W}^s_{2n}$ $(n \geq 0)$. Fix $K \geq 1$ and $\varepsilon > 0$. It is not hard to see that $\lim_{N\to\infty} \inf_{|A|\geq N} \mathbb{P}[W_n^A \neq \varnothing \ \forall n \geq 0] = 1$. Therefore, by Lemma 23, we can choose $s_0 \geq 0$ such that, for all $s \geq s_0$,

(4.45) $$\mathbb{P}[W_n^s \neq \varnothing \ \forall n \geq 0] \geq \mathbb{P}[Y_s \neq \underline{0}] - \tfrac{1}{2}\varepsilon.$$

Since the process $Y$ cannot move infinitely far in finite time, it is not hard to see from (4.42) that for each $s \in [s_0, s_0 + 2nT)$, we can choose $n_0 \geq 0$ such that, for all $n \geq n_0$,

(4.46) $$\mathbb{P}[V_{s,n} \cap [-K, K] \subset \eta(Y_{s+2nT})] \geq \mathbb{P}[W_n^s \neq \varnothing \ \forall n \geq 0] - \tfrac{1}{2}\varepsilon.$$

It follows that for all $t \geq s_0 + 2n_0 T$ we can choose $s \in [s_0, s_0 + 2nT)$ and $n \geq n_0$ such that $s + 2nT = t$ and

(4.47) $$\mathbb{P}[V_{s,n} \cap [-K, K] \subset \eta(Y_t)] \geq \mathbb{P}[Y_t \neq \underline{0} \ \forall t \geq 0] - \varepsilon.$$

Since $K$ and $\varepsilon$ are arbitrary, for each $t_k \to \infty$, we can find $K_k \to \infty$, $\varepsilon_k \to 0$, and $s_k, n_k$ such that

(4.48) $$\mathbb{P}[V_{s_k, n_k} \cap [-K_k, K_k] \subset \eta(Y_{t_k})] \geq \mathbb{P}[Y_t \neq \underline{0} \ \forall t \geq 0] - \varepsilon_k,$$

which proves our claim. $\square$

PROOF OF LEMMA 21. Set $\rho := \mathbb{P}[Y_t \neq \underline{0} \ \forall t \geq 0]$ and $\rho' := \mathbb{P}[Y_t' \neq \underline{0} \ \forall t \geq 0]$. It suffices to prove the convergence in (4.28) along an arbitrary sequence of times $t_k \to \infty$. It is clear that the limsup is bounded from above by $\rho \rho'$. To bound the liminf from below, by Lemma 22, it suffices to prove that

(4.49) $$\liminf_{k \to \infty} \mathbb{P}[D_{2L}(Y_{t_k - 1}, Y_{t_k}') \geq N] \geq \rho \rho' \qquad (N \geq 1).$$

By Lemma 24, we can couple $(Y_{t_k - 1}, Y_{t_k}')$ to random variables $V_k, V_k'$ each having law $\overline{\nu}_W$, in such a way that $(Y_{t_k - 1}, V_k)$ is independent of $(Y_{t_k}', V_k')$, and for $K \geq 1$,

(4.50) $$\liminf_{k \to \infty} \mathbb{P}[V_k \cap [-K, K] \subset \eta(Y_{t_k - 1}), V_k' \cap [-K, K] \subset \eta(Y_{t_k})] \geq \rho \rho'.$$

It is not hard to see that if $V, V'$ is a pair of independent random variables each having law $\nu_W$, then $|V \cap V'| = \infty$ a.s. Therefore, for each $N, \varepsilon > 0$, we can choose $K$ large enough such that

(4.51) $$\mathbb{P}[|V \cap V' \cap [-K, K]| < N] \leq \varepsilon,$$

and hence, by (4.50),

(4.52) $$\liminf_{k \to \infty} \mathbb{P}[D_{2L}(Y_{t_k - 1}, Y_{t_k}') \geq N] \geq \rho \rho' - \varepsilon.$$

Since $\varepsilon > 0$ is arbitrary, (4.49) follows. $\square$

**Acknowledgments.** J. Swart thanks Alison Etheridge and the Department of Statistics of Oxford University for their hospitality during a visit in October 2004. A. Sturm thanks the mathematical institutes at the universities of Erlangen-Nürnberg and Tübingen for their hospitality. We thank Maury Bramson, Ted Cox, Rick Durrett, Nina Gantert, Wolfgang König, Tom Mountford and Italo Simonelli for answering questions about their work, once more Rick Durrett for drawing our attention to reference [14], and the referee for doing a fast and thorough job.

DEPARTMENT OF MATHEMATICAL SCIENCES
UNIVERSITY OF DELAWARE
501 EWING HALL
NEWARK, DELAWARE 19716-2553
USA
E-MAIL: sturm@math.udel.edu

ÚTIA
POD VODÁRENSKOU VĚŽÍ
18208 PRAHA 8
CZECH REPUBLIC
E-MAIL: swart@utia.cas.cz